\documentclass[a4paper, oneside]{article}
\usepackage[english]{babel}
\usepackage[usenames,dvipsnames]{color}
\usepackage{tikz, amsmath, amscd, amssymb, amsfonts, stmaryrd, vmargin, amsthm, soul, bm, graphicx, float, exscale, relsize, etoolbox, mathrsfs, url}
\usepackage{hyperref}
\hypersetup{
    colorlinks,
    citecolor=blue,
    filecolor=blue,
    linkcolor=blue,
    urlcolor=blue
}
\patchcmd{\chapter}{\if@openright\cleardoublepage\else\clearpage\fi}{}{}{}

\theoremstyle{plain}        \newtheorem{theo}{Theorem}

\theoremstyle{plain}        \newtheorem{prop}[theo]{Proposition}

\theoremstyle{plain}        \newtheorem{lemm}[theo]{Lemma}

\theoremstyle{plain}        \newtheorem{coro}[theo]{Corollary}

\theoremstyle{plain}        \newtheorem{rema}[theo]{Remark}

\theoremstyle{plain}        

\theoremstyle{plain}

\begin{document}
\title{Holomorphic isometries from the Poincar\'e disk into bounded symmetric domains of rank at least two}

\author{Shan Tai Chan and Yuan Yuan\footnote{The second author is supported by National Science Foundation grant DMS-1412384, Simons Foundation grant (\#429722 Yuan Yuan) and CUSE Grant Program at Syracuse University.}}
\date{ }
\maketitle
\begin{abstract}
We first study holomorphic isometries from the Poincar\'e disk into the product of the unit disk and the complex unit $n$-ball for $n\ge 2$.
On the other hand, we observe that there exists a holomorphic isometry from the product of the unit disk and the complex unit $n$-ball into any irreducible bounded symmetric domain of rank $\ge 2$ which is not biholomorphic to any type-$\mathrm{IV}$ domain.
In particular, our study provides many new examples of holomorphic isometries from the Poincar\'e disk into irreducible bounded symmetric domains of rank at least $2$ except for type-$\mathrm{IV}$ domains.
\end{abstract}

\section{Introduction}
Holomorphic isometries between bounded symmetric domains have been studied extensively since the fundamental works of Calabi \cite{Ca53}, Clozel-Ullmo \cite{CU03} and Mok \cite{Mok12} and various rigidity results were derived. For example, when the source is irreducible and of rank at least $2$, the total geodesy of holomorphic isometries follows from the proof of Mok's theorem of metric rigidity (cf. \cite{CU03, Mok12});  when the source is of rank $1$ and of complex dimension at least $2$, namely, the complex unit $n$-ball for $n\ge 2$, and the target is the product of complex unit balls, the total geodesy of holomorphic isometries is obtained by Zhang and the second author \cite{YZ12}. On the other hand, the non-standard (i.e., not totally geodesic) holomorphic isometries were discovered by Mok either from the (complex) unit disk to polydisks \cite{Mok12} or from the complex unit $m$-ball to irreducible bounded symmetric domains of rank at least $2$ for some integer $m\ge 2$ \cite{Mok16}. It turns out to be a highly non-trivial problem to classify holomorphic isometries from the (complex) unit disk to polydisks with respect to the canonical K\"ahler-Einstein metrics (cf. \cite{Mok09, Ng10, Ch16, Ch17}). One motivation of our study is along this line after \cite{Ng10} and \cite{YZ12} to understand the holomorphic isometries from the (complex) unit disk to the product of complex unit balls with respect to the canonical K\"ahler-Einstein metrics. The first main result in the present article is that we fully characterize all holomorphic isometries from the unit disk to the product of the unit disk and a complex unit $n$-ball for $n \geq 2$ (cf. Theorem \ref{Thm_Reduction1} and Theorem \ref{Thm:Red_Para1}).
In particular, we further construct a real $1$-parameter family of mutually incongruent holomorphic isometries.
When the target is an irreducible bounded symmetric domain $\Omega$ of rank at least $2$, holomorphic isometries from the unit disk are poorly understood in general.
This is the second motivation of our study as one may obtain various examples by composing the maps described above with a totally geodesic (or non-standard) holomorphic isometric embedding from the product of the unit disk with a complex unit $n$-ball to $\Omega$.
When the target $\Omega$ is of rank at least $3$, by composing the real $1$-parameter family of mutually incongruent holomorphic isometries from the unit disk to the $3$-disk found by Mok (cf. \cite{Mok12}) with the totally geodesic holomorphic isometric embedding from the $3$-disk to $\Omega$, it is known that a real $1$-parameter family of mutually incongruent holomorphic isometries from the unit disk to $\Omega$ exists.
When $\mathrm{rank}(\Omega)=2$, the real $1$-parameter family of mutually incongruent holomorphic isometries from the unit disk to $\Omega$ is only known to exist when $\Omega$ is a classical domain of type $\mathrm{IV}$ (cf. \cite{CM17, XY16b}). Nevertheless, we may obtain a real $1$-parameter family of mutually incongruent holomorphic isometries from the unit disk to any classical irreducible bounded symmetric domain $\Omega$ of rank at least $2$ except for type-$\mathrm{IV}$ domains (cf. Theorem \ref{pro:FamHI_MuIC1_ClassicalDomain}). Another interesting question of Mok regards the boundary extension of holomorphic isometries from the unit disk to bounded symmetric domains of rank at least $2$ (cf. Problem 5.2.2 in \cite{Mok11}).
In this direction, we construct new non-standard holomorphic isometries from the unit disk into $\Omega$ that extend holomorphically to a neighborhood of the closed unit disk and have irrational component function(s), where $\Omega$ is any irreducible bounded symmetric domain of rank at least $2$ in its Harish-Chandra realization (cf. Theorem \ref{Pro:NewEx_ExHolo_irrat}). Note that rational examples are known to exist before (cf. \cite{CM17, XY16b}). For related problems of holomorphic isometries between Hermitian symmetric spaces of compact type, the interested readers may refer to \cite{HY14,Eb15,Y16}.

\bigskip
We fix the notations in the present article.
Denote by $g_D$ the canonical K\"ahler-Einstein metric on an irreducible bounded symmetric domain $D\Subset \mathbb C^n$ normalized so that minimal disks are of constant Gaussian curvature $-2$.
We also denote by $ds_U^2$ the Bergman metric of any bounded domain $U\Subset \mathbb C^N$. On an irreducible bounded symmetric domain $D\Subset \mathbb C^n$, $g_D$ and $ds_D^2$ are proportional by a positive constant. In particular, on a complex unit $n$-ball $\mathbb{B}^n$, $ds_{\mathbb{B}^n}^2 = (n+1) g_{\mathbb{B}^n}$.
Let $D\Subset \mathbb C^n$ and $\Omega\Subset \mathbb C^N$ be bounded symmetric domains.
We write $D=D_1\times\cdots \times D_m$ and $\Omega=\Omega_1\times\cdots \times \Omega_k$, where $D_j$, $1\le j\le m$, (resp. $\Omega_l$, $1\le l\le k$), are irreducible factors of $D$ (resp. $\Omega$).
Let $F_1=(f^{(1)}_1, \cdots, f^{(1)}_k),F_2=(f^{(2)}_1, \cdots, f^{(2)}_k):(D_1,\lambda_1g_{D_1})\times\cdots\times (D_m,\lambda_m g_{D_m})
\to (\Omega_1,\lambda'_1g_{\Omega_1})\times\cdots \times (\Omega_k,\lambda'_kg_{\Omega_k})$ be holomorphic isometries for some positive real constants $\lambda_j$, $1\le j\le m$, and $\lambda'_j$, $1\le j\le k$,
in the sense
\[ \bigoplus_{j=1}^m \lambda_j g_{D_j} = \sum_{j=1}^k \lambda'_j (f^{(l)}_j)^* g_{\Omega_j}\]
for $l=1,2$.
Then, $F_1$ and $F_2$ are said to be congruent (or equivalent) to each other if $F_1= \Psi\circ F_2 \circ \phi$ for some $\phi\in \mathrm{Aut}(D)$ and $\Psi\in \mathrm{Aut}(\Omega)$; otherwise, $F_1$ and $F_2$ are said to be incongruent to each other.
This defines an equivalence class of any holomorphic isometry $F:(D_1,\lambda_1g_{D_1})\times\cdots\times (D_m,\lambda_m g_{D_m})
\to (\Omega_1,\lambda'_1g_{\Omega_1})\times\cdots \times (\Omega_k,\lambda'_kg_{\Omega_k})$.

Throughout the present article, we say that $b\in \mathbb P^1$ is a branch point of a rational function $R:\mathbb P^1\to \mathbb P^1$ of degree $\deg(R)\ge 2$ if $b=R(a)$ for some ramification point $a$ of $R$, where $R:\mathbb P^1\to \mathbb P^1$ is regarded as a finite branched covering of degree equal to $\deg(R)\ge 2$.

\bigskip
 
{\bf Acknowledgement} We would like to thank Prof. X. Huang and Prof. N. Mok for their interest and invaluable comments. We are grateful to the referee for carefully reading the manuscript and the helpful suggestions. 

\section{Holomorphic isometries from the unit disk to $\Delta\times \mathbb B^n$ for $n\ge 2$}
For any integer $n\ge 1$, we denote by $\mathbb B^n$ the complex unit $n$-ball in the complex $n$-dimensional Euclidean space $\mathbb C^n$, i.e., 
\[ \mathbb B^n:=\left\{(z_1,\ldots,z_n)\in \mathbb C^n:\sum_{j=1}^n |z_j|^2 <1\right\}. \]
Then, the K\"ahler form $\omega_{g_{\mathbb B^n}}$ of $(\mathbb B^n,g_{\mathbb B^n})$, $n\ge 1$, is given by
\[ \omega_{g_{\mathbb B^n}}=-\sqrt{-1}\partial\overline\partial \log \left(1-\sum_{j=1}^n |z_j|^2\right). \]
We denote by $\Delta:=\mathbb B^1\Subset \mathbb C$ the (open) unit disk in the complex plane $\mathbb C$ throughout the present article.
For $n\ge 2$, we define the complex-analytic subvariety
\[ W_{{\bf U'}}:= \left\{(w,z_1,\ldots,z_n)\in \Delta\times \mathbb B^n:{\bf U'} \begin{pmatrix}
w,z_1,\ldots, z_n
\end{pmatrix}^T
= \begin{pmatrix}
wz_1,\ldots, wz_n
\end{pmatrix}^T \right\} \]
of $\Delta\times \mathbb B^n$ for any matrix ${\bf U'} \in M(n,n+1;\mathbb C)$ of full rank $n$.

The following extension theorem is a special case of Calabi's theorem (cf. \cite{Ca53}).
\begin{prop}\label{Pro_ExThm}
Let $f:(\Delta,g_\Delta;0)\to (\Delta,g_\Delta;0)\times (\mathbb B^n,g_{\mathbb B^n};{\bf 0})$ be a germ of holomorphic isometry, where $n\ge 2$ is an integer.
Then, $f$ extends to a holomorphic isometric embedding
$F:(\Delta,g_\Delta)\to (\Delta,g_\Delta)\times (\mathbb B^n,g_{\mathbb B^n})$.
\end{prop}
\begin{proof}
We may suppose that $f$ is defined on $B^1(0,\varepsilon)$ for some $\varepsilon>0$.
Writing $f=(f_1,f_{2,1},\ldots,f_{2,n})$, we have the functional equation
\[ \left(1-|f_1(w)|^2\right)\left(1-\sum_{j=1}^n |f_{2,j}(w)|^2\right) = 1-|w|^2 \]
for any $w\in B^1(0,\varepsilon)$.
In particular, we have
\[ \left(1-|f_1(w)|^2\right)^{2(n+1)}\left(1-\sum_{j=1}^n |f_{2,j}(w)|^2\right)^{2(n+1)} = (1-|w|^2)^{2(n+1)} \]
for any $w\in B^1(0,\varepsilon)$.
Write $\Omega:=\Delta^{n+1} \times \mathbb B^n \times \mathbb B^n$.
Define a holomorphic isometry $\widetilde f:(\Delta,(n+1)ds_\Delta^2)\times (\mathbb B^n,2ds_{\mathbb B^n}^2)\to (\Omega,ds_\Omega^2)$ by
$\widetilde f(w,z_1,\ldots,z_n)
= (w,\ldots,w;z_1,\ldots,z_n;z_1,\ldots,z_n)$.
Then, we see that $\widetilde f\circ f: (\Delta,(n+1)ds_\Delta^2)\to (\Omega,ds_\Omega^2)$ is a germ of holomorphic isometry.
Therefore, $\widetilde f\circ f$ extends to the holomorphic isometry $\widetilde F: (\Delta,(n+1)ds_\Delta^2)\to (\Omega,ds_\Omega^2)$ by the extension theorem for local holomorphic isometries between bounded symmetric domains (cf. Calabi \cite{Ca53} and Mok \cite{Mok12}).
In particular, $f$ extends to a holomorphic isometric embedding
$F:(\Delta,g_\Delta)\to (\Delta,g_\Delta)\times (\mathbb B^n,g_{\mathbb B^n})$.
\end{proof}
It is natural to ask whether all holomorphic isometries $(\Delta,g_\Delta)\to (\Delta,g_\Delta)\times (\mathbb B^n,g_{\mathbb B^n})$ are obtained from the square-root embedding $\Delta\to \Delta^2$ if the components $\Delta\to \Delta$ and $\Delta\to \mathbb B^n$ are not constant maps. But we will show that there are examples of holomorphic isometries which are not obtained in that way.
\subsection{Existence of holomorphic isometries}
\text{}\\
Let ${\bf U}:=\begin{pmatrix}
u_{1,1} &\cdots & u_{1,n+1}\\
\vdots &\ddots & \vdots\\
u_{n+1,1} &\cdots & u_{n+1,n+1}\\
\end{pmatrix}\in U(n+1)$ be a unitary matrix, where $n\ge 2$ is an integer.
Our goal in this section is to obtain a holomorphic isometry $f:=(f_1,f_{2,1},\ldots,f_{2,n}):(\Delta,g_\Delta)\to (\Delta,g_\Delta)\times (\mathbb B^n,g_{\mathbb B^n})$ such that
\begin{equation}\label{Eq:DefEq1}
 {\bf U} \begin{pmatrix}
f_1(w),f_{2,1}(w),\ldots, f_{2,n}(w)
\end{pmatrix}^T
= \begin{pmatrix}
w,f_1(w)  f_{2,1}(w),\ldots, f_1(w) f_{2,n}(w)
\end{pmatrix}^T
\end{equation}
for any $w\in \Delta$.
In other words, we need to solve the system provided in Eq. (\ref{Eq:DefEq1})
for some germ of holomorphic function $f_1:(\Delta;0)\to (\Delta;0)$ and some germ of holomorphic map $(f_{2,1},\ldots,f_{2,n}):(\Delta;0) \to (\mathbb B^n;{\bf 0})$.
Then, we have a germ of holomorphic isometry $f:=(f_1,f_{2,1},\ldots,f_{2,n}):(\Delta,g_\Delta;0)\to(\Delta,g_\Delta;0)\times (\mathbb B^n,g_{\mathbb B^n};{\bf 0})$ and the rest would follow from Proposition \ref{Pro_ExThm}.

Write ${\bf U'}:=\begin{pmatrix}
u_{2,1} &\cdots & u_{2,n+1}\\
\vdots &\ddots & \vdots\\
u_{n+1,1} &\cdots & u_{n+1,n+1}\\
\end{pmatrix}$.
Then, it is obvious that ${\bf U'}$ is of full rank $n$.
Since $W_{{\bf U'}}$ is a complex-analytic subvariety of $\Delta\times \mathbb B^n$ which is smooth and of dimension $1$ at ${\bf 0}$, there exists a germ of holomorphic map $f:=(f_1,f_{2,1},\ldots,f_{2,n}):(\Delta,0) \to (\Delta\times \mathbb B^n;{\bf 0})$ such that the image of $f$ is an open neighborhood of ${\bf 0}$ in $W_{\bf U'}$ and
\[ {\bf U'}
\begin{pmatrix}
f_1(w),
f_{2,1}(w),\ldots,f_{2,n}(w)
\end{pmatrix}^T
= \begin{pmatrix}
f_1(w)f_{2,1}(w),\ldots,f_1(w)f_{2,n}(w)
\end{pmatrix}^T, \]
equivalently
\begin{equation}\label{EqSystem2}
\left(\begin{pmatrix}
u_{2,2} &\cdots & u_{2,n+1}\\
\vdots &\ddots & \vdots\\
u_{n+1,2} &\cdots & u_{n+1,n+1}\\
\end{pmatrix} - f_1(w) {\bf I_n}\right)
\begin{pmatrix}
f_{2,1}(w)\\
\vdots\\
f_{2,n}(w)
\end{pmatrix}
= -f_1(w)\begin{pmatrix}
u_{2,1} \\
\vdots \\
u_{n+1,1}
\end{pmatrix}
\end{equation}
for any $w$ lying inside the domain of $f$.
In order to apply the Cramer's rule to solve the system in Eq. (\ref{EqSystem2}), we need
\[ \det\left(\begin{pmatrix}
u_{2,2}&\cdots & u_{2,n+1}\\
\vdots & \ddots & \vdots\\
u_{n+1,2}&\cdots & u_{n+1,n+1}
\end{pmatrix}-f_1(w) {\bf I_n}\right) \neq 0 \]
for any $w$ in some open neighborhood of $0$ in the domain of $f$.
Thus, it suffices to require that
\[ \det \begin{pmatrix}u_{2,2}&\cdots & u_{2,n+1}\\\vdots & \ddots & \vdots\\u_{n+1,2}&\cdots & u_{n+1,n+1}\end{pmatrix} \neq 0. \]
Actually, we have the following:
\begin{prop}
\label{Pro_Det_Constf1}
Write ${\bf U}'':=\begin{pmatrix}
u_{2,2}&\cdots & u_{2,n+1}\\
\vdots & \ddots & \vdots\\
u_{n+1,2}&\cdots & u_{n+1,n+1}
\end{pmatrix}$.
In the above settings, we have $f_1(w)\equiv 0$ if and only if $\det {\bf U''}=0$.
\end{prop}
\begin{proof}
By performing Gaussian elimination, if
$\det {\bf U''}=0$,
then we have 
\[ c f_1(w) = f_1(w) \sum_{j=1}^n c_j f_{2,j}(w) \]
for some $c,c_j\in \mathbb C$, $1\le j\le n$, such that $c\neq 0$.
Then, we have
\[ f_1(w) \left( c- \sum_{j=1}^n c_j f_{2,j}(w) \right) = 0, \]
so that either $f_1(w)\equiv 0$ or $c- \sum_{j=1}^n c_j f_{2,j}(w)\equiv 0$. But the latter is impossible because $c\neq 0$ and $f_{2,j}(0)=0$. Thus, $f_1(w)\equiv 0$.

Conversely, if $f_1(w)\equiv 0$, then
${\bf U''}
\begin{pmatrix}
f_{2,1}(w),\ldots,f_{2,n}(w)
\end{pmatrix}^T
\equiv {\bf 0}$ by Eq. (\ref{EqSystem2}).
Since there is $j$, $1\le j\le n$, such that $f_{2,j}(w)\not\equiv 0$, the matrix ${\bf U''}$ is not invertible, i.e., $\det {\bf U''} =0$.
In this case, we can always solve
\begin{equation}\label{determinant}
 {\bf U}
\begin{pmatrix}
0 ,
f_{2,1}(w),\ldots,
f_{2,n}(w)
\end{pmatrix}^T
= \begin{pmatrix}
w,0,\ldots,0
\end{pmatrix}^T.
\end{equation}
Actually, we have 
$\begin{pmatrix}
0 ,
f_{2,1}(w),\ldots,
f_{2,n}(w)
\end{pmatrix}^T= {\bf U}^{-1} \begin{pmatrix}
w,0,\ldots,0
\end{pmatrix}^T$.
The proof is complete.
\end{proof}
\noindent We are ready to prove the following existence theorem for any unitary matrix ${\bf U}:=\begin{pmatrix}
u_{ij}
\end{pmatrix}_{1\le i,j\le n+1}$ $\in$ $U(n+1)$.
\begin{theo}[Existence Theorem]\label{Thm_Existence1}
Let ${\bf U}:=\begin{pmatrix}
u_{ij}
\end{pmatrix}_{1\le i,j\le n+1} \in U(n+1)$ be a unitary matrix, where $n\ge 2$ is an integer.
Then, there is a holomorphic isometry
$f:=(f_1,f_{2,1},\ldots,f_{2,n}):(\Delta,g_\Delta)\to (\Delta,g_\Delta)\times (\mathbb B^n,g_{\mathbb B^n})$ such that
\[  {\bf U} \begin{pmatrix}
f_1(w),f_{2,1}(w),\ldots, f_{2,n}(w)
\end{pmatrix}^T
= \begin{pmatrix}
w,f_1(w)  f_{2,1}(w),\ldots, f_1(w) f_{2,n}(w)
\end{pmatrix}^T \]
for any $w\in \Delta$.
\end{theo}
\begin{proof}
We adapt the above settings and constructions. In particular, there exists a germ of holomorphic map $f:=(f_1,f_{2,1},\ldots,f_{2,n}):(\Delta;0) \to (\Delta\times \mathbb B^n;{\bf 0})$ such that
\begin{equation}\label{Eq:System1}
\left(\begin{pmatrix}
u_{2,2} &\cdots & u_{2,n+1}\\
\vdots &\ddots & \vdots\\
u_{n+1,2} &\cdots & u_{n+1,n+1}\\
\end{pmatrix} - f_1(w) {\bf I_n}\right)
\begin{pmatrix}
f_{2,1}(w)\\
\vdots\\
f_{2,n}(w)
\end{pmatrix}
= -f_1(w)\begin{pmatrix}
u_{2,1} \\
\vdots \\
u_{n+1,1}
\end{pmatrix}
\end{equation}
on the domain of $f$.
Write ${\bf U''}:=\begin{pmatrix}
u_{2,2} &\cdots & u_{2,n+1}\\
\vdots &\ddots & \vdots\\
u_{n+1,2} &\cdots & u_{n+1,n+1}\\
\end{pmatrix}$.
If $\det {\bf U''}=0$, then by Proposition \ref{Pro_Det_Constf1} we may solve $f(w)=(0,f_{2,1}(w),\ldots,f_{2,n}(w))$ out from Eq. (\ref{determinant}) and the result follows.

From now on, we assume that $\det {\bf U''} \neq 0$.
Then, by applying the Cramer's rule to Eq. (\ref{Eq:System1}), there exists a rational function $R_j:\mathbb P^1\to \mathbb P^1$ such that
$f_{2,j}(w)=R_j(f_1(w))$ on $B^1(0,\varepsilon')$ for some $\varepsilon'>0$ and any pole of $R_j$ in $\mathbb C$ is a zero of $\det \left( {\bf U}'' - z {\bf I_n}\right)$ for $1\le j\le n$.

Now, we require that the germ of holomorphic map $f:=(f_1,f_{2,1},\ldots,f_{2,n}):(\Delta,0) \to (\Delta\times \mathbb B^n;{\bf 0})$ satisfies
\begin{equation}\label{Eq:ExtraEqforgerms}
u_{11} f_1(w) + \sum_{j=1}^n u_{1,j+1} f_{2,j}(w) = w.
\end{equation}
This is equivalent to the requirement that the germ $f_1:\Delta\to \Delta$ of holomorphic function satisfies $R(f_1(w))=w$,
where $R:\mathbb P^1\to \mathbb P^1$ is the rational function defined by
\begin{equation}\label{Eq:RatFun1}
R(z) = {z \det \left({\bf U} - z \begin{bmatrix}
{\bf 0} & {\bf 0}\\
{\bf 0} & {\bf I_n}
\end{bmatrix}\right)\over \det \left( {\bf U''} - z {\bf I_n}\right)}. 
\end{equation}
This follows from the fact that $f_{2,j}=R_j\circ f_1$ on $B^1(0,\varepsilon')$ for some rational function $R_j:\mathbb P^1\to \mathbb P^1$, $1\le j\le n$,
and by substituting $f_{2,j}(w)=R_j(f_1(w))$, $1\le j\le n$, into Eq. (\ref{Eq:ExtraEqforgerms}).
Since $\det{\bf U}\neq 0$ and $\det {\bf U''}\neq 0$, $z=0$ is a simple zero of $R$ so that $0$ is not a ramification point of $R$.
In particular, there are neighborhoods $U$ and $V$ of $0$ in $\mathbb C$ such that $R|_U:U\to V$ is a biholomorphism.
Therefore, there is a germ of holomorphic function $f_1:(\Delta;0)\to (\Delta;0)$ such that $R(f_1(w))=w$ on the domain of $f_1$.
Then, we have
\[  {\bf U} \begin{pmatrix}
f_1(w),f_{2,1}(w),\ldots, f_{2,n}(w)
\end{pmatrix}^T
= \begin{pmatrix}
w,f_1(w)  f_{2,1}(w),\ldots, f_1(w) f_{2,n}(w)
\end{pmatrix}^T \]
for some germ of holomorphic function $f_1:(\Delta;0)\to (\Delta;0)$ and some germ of holomorphic map $(f_{2,1},\ldots,f_{2,n}):(\Delta;0) \to \mathbb B^n$ satisfying
$f_{2,j}=R_j\circ f_1$ for $1\le j\le n$.
It follows that
\[ (1-|f_1(w)|^2)\left(1-\sum_{j=1}^n |f_{2,j}(w)|^2\right)
= 1-|w|^2 \]
on $B^1(0,\varepsilon'')$ for some $\varepsilon''>0$.
Thus, $f:=(f_1,f_{2,1},\ldots,f_{2,n}):(\Delta,g_\Delta;0)\to (\Delta,g_\Delta;0)\times(\mathbb B^n,g_{\mathbb B^n};{\bf 0})$ is a germ of holomorphic isometry.
It follows from Proposition \ref{Pro_ExThm} that $f$ extends to a holomorphic isometry
$(\Delta,g_\Delta)\to (\Delta,g_\Delta)\times (\mathbb B^n,g_{\mathbb B^n})$, which is also denoted by $f$.
The rest follows directly from the above constructions.
\end{proof}
\begin{rema}\label{remark2.4}
\begin{enumerate}
\item[(1)]
This theorem applies to the case where $n=1$ and this would reduce to the case of holomorphic isometries $(\Delta,ds_\Delta^2)\to (\Delta^2,ds_{\Delta^2}^2)$.
\item[(2)] 
Writing ${\bf U}''=\begin{pmatrix}u_{2,2}&\cdots & u_{2,n+1}\\\vdots & \ddots & \vdots\\u_{n+1,2}&\cdots & u_{n+1,n+1}\end{pmatrix}$, we also have
\[ R(z) = {z u_{11} \det\left( {\bf U}'' - {1\over u_{11}}\begin{pmatrix}u_{21}\\\vdots \\ u_{n+1,1}\end{pmatrix}\cdot \begin{pmatrix}u_{12},\ldots, u_{1,n+1} \end{pmatrix} - z {\bf I_n} \right)\over \det \left( {\bf U}'' - z {\bf I_n}\right)}, \]
when ${\bf U}''$ is invertible.
\item[(3)] Let ${\bf U}:=\begin{pmatrix}
u_{ij}
\end{pmatrix}_{1\le i,j\le n+1} \in U(n+1)$ be any unitary matrix.
If the matrix ${\bf U}''=\begin{pmatrix}u_{ij}\end{pmatrix}_{2\le i,j\le n+1}$ is invertible, then one constructs the rational function $R:\mathbb P^1\to \mathbb P^1$ from ${\bf U}$ so that $R$ determines the component function $f_1:\Delta\to\Delta$ of the holomorphic isometry $f:=(f_1,f_{2,1},\ldots,f_{2,n}):(\Delta,g_\Delta)\to (\Delta,g_\Delta)\times (\mathbb B^n,g_{\mathbb B^n})$ uniquely via the identity $R(f_1(w))=w$.
This follows from the fact that $R$ is unramified at $0$ and this gives rise to a unique germ of holomorphic function $f_1:(\Delta;0)\to (\Delta;0)$ such that $R(f_1(w))=w$.
From the functional equation, it is obvious that if there is another holomorphic isometry $\widetilde f:=(f_1,g_1,\ldots,g_n):(\Delta,g_\Delta)\to (\Delta,g_\Delta)\times (\mathbb B^n,g_{\mathbb B^n})$ such that $\widetilde f(0)={\bf 0}$, then $\sum_{j=1}^n |g_j(w)|^2 = \sum_{j=1}^n |f_{2,j}(w)|^2$ for any $w\in \Delta$ so that $\widetilde f$ and $f$ are congruent to each other. It turns out that the rational function $R:\mathbb P^1\to \mathbb P^1$ determines the holomorphic isometry $f:=(f_1,f_{2,1},\ldots,f_{2,n}):(\Delta,g_\Delta)\to (\Delta,g_\Delta)\times (\mathbb B^n,g_{\mathbb B^n})$ up to congruence.

Conversely, let $f:=(f_1,f_{2,1},\ldots,f_{2,n}):(\Delta,g_\Delta)$ $\to$ $(\Delta,g_\Delta)\times (\mathbb B^n,g_{\mathbb B^n})$ be any holomorphic isometry. We may assume that $f(0)={\bf 0}$.
Then, it follows from the functional equation and Calabi's theorem (cf. \cite[Theorem 2]{Ca53}) that
$${\bf U} \begin{pmatrix}
f_1(w),f_{2,1}(w),\ldots, f_{2,n}(w)
\end{pmatrix}^T
= \begin{pmatrix}
w,f_1(w)  f_{2,1}(w),\ldots, f_1(w) f_{2,n}(w)
\end{pmatrix}^T$$
for some unitary matrix ${\bf U}\in U(n+1)$.
If $f_1$ is non-constant, then we can construct a unique rational function $R:\mathbb P^1\to \mathbb P^1$ from ${\bf U}$ so that $R(f_1(w))=w$.
\end{enumerate}
\end{rema}
\subsubsection{Normalization of matrices}\label{Sec:2.1.1}
Let $f=(f_1,f_{2,1},\ldots,f_{2,n}):(\Delta,g_\Delta)\to (\Delta,g_\Delta)\times(\mathbb B^n,g_{\mathbb B^n})$ be a holomorphic isometry such that $f_1$ is a non-constant function, where $n\ge 2$ is an integer.
Assume without loss of generality that $f(0)={\bf 0}$.
Then, we have
\[ \left(1-|f_1(w)|^2\right)\left(1-\sum_{j=1}^n |f_{2,j}(w)|^2\right) = 1-|w|^2 \]
for any $w\in \Delta$.
It follows from the local rigidity theorem of Calabi \cite[Theorem 2]{Ca53} that there is ${\bf U}\in U(n+1)$ such that
\[  {\bf U} \begin{pmatrix}
f_1(w),f_{2,1}(w),\ldots, f_{2,n}(w)
\end{pmatrix}^T
= \begin{pmatrix}
w,f_1(w)  f_{2,1}(w),\ldots, f_1(w) f_{2,n}(w)
\end{pmatrix}^T \]
for any $w\in \Delta$.
Let ${\bf B}\in U(n)$ be a unitary matrix.
Define a holomorphic map $(g_1,\ldots,g_n):\Delta\to \mathbb B^n$ by
\[ \begin{pmatrix}
g_1(w),\ldots,g_n(w)
\end{pmatrix}^T={\bf B}\begin{pmatrix}
f_{2,1}(w),\ldots, f_{2,n}(w)
\end{pmatrix}^T. \]
Then, $(f_1,g_1,\ldots,g_n):(\Delta,g_\Delta)\to (\Delta,g_\Delta)\times(\mathbb B^n,g_{\mathbb B^n})$ is a holomorphic isometry which is congruent to the holomorphic isometry $f:(\Delta,g_\Delta)\to (\Delta,g_\Delta)\times(\mathbb B^n,g_{\mathbb B^n})$.
Moreover, we have
\[ \begin{bmatrix}
1 & \\
 & {\bf B}
\end{bmatrix}
{\bf U} \begin{bmatrix}
1 & \\
 & {\bf B}^{-1}
\end{bmatrix}  \begin{pmatrix}
f_1(w)\\
g_1(w)\\
\vdots\\
g_n(w)
\end{pmatrix}
= \begin{pmatrix}
w\\
f_1(w)g_1(w)\\
\vdots\\
f_1(w)g_n(w)
\end{pmatrix}. \]
Here we choose ${\bf B}\in U(n)$ so that
$${\bf B} \begin{pmatrix}
u_{22}&\cdots & u_{2,n+1}\\
\vdots &\ddots &\vdots\\
u_{n+1,2}&\cdots & u_{n+1,n+1}
\end{pmatrix} {\bf B}^{-1}$$ is an upper triangular matrix by the Schur Decomposition (cf. Theorem 3.3 in \cite[p.\;79]{Zh11}).
(Noting that ${\bf B}^{-1}=\overline{{\bf B}}^T$.)
In particular, we may write
\[ \begin{bmatrix}
1 & \\
 & {\bf B}
\end{bmatrix}
{\bf U} \begin{bmatrix}
1 & \\
 & {\bf B}^{-1}
\end{bmatrix}
= \begin{pmatrix}
u_{11}' & u_{12}' & \cdots & u_{1,n+1}'\\
u_{21}' & u_{22}' &\cdots & u_{2,n+1}'\\
\vdots & \vdots & \ddots&\vdots\\
u_{n+1,1}' & 0 & \cdots & u_{n+1,n+1}'
\end{pmatrix}, 
\]
where $u_{11}'=u_{11}$, $\begin{pmatrix}
u_{12}' ,\ldots, u_{1,n+1}'
\end{pmatrix} = \begin{pmatrix}
u_{12},\ldots, u_{1,n+1}
\end{pmatrix} {\bf B}^{-1}$,
\[ \begin{pmatrix}
u_{21}' \\
\vdots\\
u_{n+1,1}'
\end{pmatrix} = {\bf B}\begin{pmatrix}
u_{21} \\
\vdots\\
u_{n+1,1}
\end{pmatrix} \]
and
\[ \begin{pmatrix}
u_{22}' &\cdots & u_{2,n+1}'\\
 & \ddots&\vdots\\
{\bf 0} &  & u_{n+1,n+1}'
\end{pmatrix} = {\bf B} \begin{pmatrix}
u_{22}&\cdots & u_{2,n+1}\\
\vdots &\ddots &\vdots\\
u_{n+1,2}&\cdots & u_{n+1,n+1}
\end{pmatrix} {\bf B}^{-1}. \]
Then, $$\begin{pmatrix}
u_{22}' &\cdots & u_{2,n+1}'\\
 & \ddots&\vdots\\
{\bf 0} &  & u_{n+1,n+1}'
\end{pmatrix}$$ is invertible if and only if $$\begin{pmatrix}
u_{22}&\cdots & u_{2,n+1}\\
\vdots &\ddots &\vdots\\
u_{n+1,2}&\cdots & u_{n+1,n+1}
\end{pmatrix}$$ is invertible.
From now on, we may use the normalization of the unitary matrix directly, i.e., $${\bf U}'':=\begin{pmatrix}
u_{2,2}&\cdots & u_{2,n+1}\\
\vdots & \ddots & \vdots\\
u_{n+1,2}&\cdots & u_{n+1,n+1}
\end{pmatrix}$$ is assumed to be upper triangular.

\subsection{Degree of the rational function $R$}
We observe that the rational function $R$ is of a certain special form.
\begin{lemm}\label{lem:RationalF1}
Let $f=(f_1,f_{2,1},\ldots,f_{2,n}): (\Delta,g_\Delta)\to (\Delta,g_\Delta)\times(\mathbb B^n,g_{\mathbb B^n})$ be a holomorphic isometry such that $f_1$ is a non-constant function and $f(0)={\bf 0}$, where $n\ge 2$.
Then, there is a rational function $R:\mathbb P^1\to \mathbb P^1$ such that $R(f_1(w))=w$, $R\left({1\over \overline z}\right)={1\over \overline{R(z)}}$ and
\[ R(z)=\alpha_0 z\prod_{j=1}^n {z-{1\over \overline{\alpha_j}}\over z-\alpha_j}, \]
where $\alpha_j\in \overline\Delta \smallsetminus \{0\}$ for $1\le j\le n$ and $\alpha_0\in \overline\Delta\smallsetminus \{0\}$.
\end{lemm}
\begin{proof}
By polarization, it follows from the functional equation that
\[ (1-f_1(w)\overline{f_1(\zeta)})\left(1-\sum_{j=1}^n f_{2,j}(w)\overline{f_{2,j}(\zeta)}\right) = 1-w\overline\zeta \]
for any $w,\zeta\in \Delta$.
It follows also from the local rigidity theorem of Calabi \cite[Theorem 2]{Ca53} that
\[ {\bf U} \begin{pmatrix}
f_1(w),f_{2,1}(w),\ldots, f_{2,n}(w)
\end{pmatrix}^T
= \begin{pmatrix}
w , f_1(w) f_{2,1}(w),\ldots, f_1(w)f_{2,n}(w)
\end{pmatrix}^T \]
for some ${\bf U}\in U(n+1)$.
By the normalization of the matrix ${\bf U}$, we can assume that ${\bf U}=\begin{bmatrix}
u_{11} & {\bf u}\\
{\bf v} & {\bf U''}
\end{bmatrix}$ with ${\bf U''}$ being an $n$-by-$n$ upper triangular matrix.
From the assumption that $f_1$ is non-constant, we have $u_{jj}\neq 0$ for $1\le j\le n+1$ by Proposition \ref{Pro_Det_Constf1}.
Then, we have $f_{2,j}(w)=R_j(f_1(w))$ for some rational function $R_j:\mathbb P^1\to \mathbb P^1$ for $1\le j\le n$.
Moreover, we have $R(f_1(w))=w$ for some rational function $R:\mathbb P^1\to \mathbb P^1$
of the form
$R(z) = \gamma z{p(z)\over \prod_{j=2}^{n+1} (z-u_{jj})}$ from the construction,
where $p(z)$ is a complex polynomial and $\gamma\in \mathbb C$ is a nonzero constant such that $p(0)\neq 0$.
From the polarized functional equation, we have
\[ (1-f_1(w)\overline{f_1(\zeta)})\left(1-\sum_{j=1}^n 
R_j(f_1(w))\overline{R_j(f_1(\zeta))}\right) = 1-R(f_1(w))\overline{R(f_1(\zeta))} \]
for any $w,\zeta\in \Delta$.
Since $f_1:\Delta\to\Delta$ is a non-constant holomorphic function, the image of $f_1$ is an open subset of $\mathbb C\subset \mathbb P^1$ by the Open Mapping Theorem.
Thus, we have 
\[ (1-\xi\overline{\eta})\left(1-\sum_{j=1}^n 
R_j(\xi)\overline{R_j(\eta)}\right) = 1-R(\xi)\overline{R(\eta)} \]
for any $\xi,\eta\in \mathbb C\smallsetminus\{u_{jj}\mid 2\le j\le n+1\}$.
In particular, we have
\begin{equation}\label{Eq:FE_R}
(1-|\xi|^2)\left(1-\sum_{j=1}^n 
|R_j(\xi)|^2\right) = 1-|R(\xi)|^2\
\end{equation}
for any $\xi\in \mathbb C\smallsetminus\{u_{jj}\mid 2\le j\le n+1\}$ (cf. \cite{Ch17}).
Then, there is an open arc $A\subset \partial\Delta$ such that $|R_j(\xi_0)|^2$, $1\le j\le n$, are finite for any $\xi_0\in A$ because the set of poles of $R_j$ in $\mathbb C$ is a subset of $\{u_{ll}:2\le l\le n+1\}$ for $1\le j\le n$.
Then, we have $|R(\xi)|^2=1$ for any $\xi\in A$ by Eq. (\ref{Eq:FE_R}).
We define the rational function $\Upsilon:\mathbb P^1\to \mathbb P^1$ by $\Upsilon(z):=R(z)\overline{R\left({1\over \overline{z}}\right)}$.
Then, we have $\Upsilon(z)= 1$ for any $z\in A$.
But then $\Upsilon^{-1}(1)$ is finite if $\Upsilon$ is non-constant.
Therefore, $\Upsilon(z)\equiv 1$ is a constant function.
In particular, we have $\overline{R\left({1\over \overline{z}}\right)} = {1\over R(z)}$.
Now, $z_0$ is a zero of $R$ if and only if ${1\over \overline{z_0}}$ is a pole of $R$.
The poles of $R$ in $\mathbb P^1=\mathbb C\cup\{\infty\}$ are precisely the infinity $\infty$ and $u_{ll}$ for $2\le l\le n+1$ such that $|u_{ll}|\neq 1$.
Thus, the zeros of $R$ in $\mathbb P^1$ are precisely $0$ and ${1\over \overline{u_{ll}}}$ for $2\le l\le n+1$ such that $|u_{ll}|\neq 1$.
In particular, one has
\[ R(z) = \gamma' z
\prod_{j=2}^{n+1}{z-{1\over \overline{u_{jj}}}\over z- u_{jj}} \]
for some $\gamma'\in \mathbb C\smallsetminus\{0\}$.
(Noting that if $|u_{jj}|=1$ for some $j$, $2\le j\le n+1$, then ${z-{1\over \overline{u_{jj}}}\over z- u_{jj}}\equiv 1$.)
Since $\overline{R\left({1\over \overline{z}}\right)} = {1\over R(z)}$, we have $|\gamma'|=\prod_{j=2}^{n+1} |u_{jj}|$.
Comparing $R(z)$ with the formula that we have obtained in item $(2)$ of Remark \ref{remark2.4}, we have $\gamma'=u_{11}$.
\end{proof}
\begin{rema}\label{remark2.6}
We have $\deg(R)\le n+1$.
Moreover, if $|u_{jj}|^2<1$ for $2\le j \le n+1$, then we have $\deg(R)=n+1$.
\end{rema}
\noindent This actually yields the following corollary.
\begin{coro}
Let ${\bf U}=\begin{pmatrix}
u_{ij}
\end{pmatrix}_{1\le i,j\le n}\in U(n)$ be such that
$u_{kj}=0$ for $2\le j\le k-1$, $3\le k\le n$ and $u_{11}\neq 0$, where $n\ge 3$ is an integer.
Then, the zeros of the polynomial
\[ p_{\bf U}(z):= \det \left(
{\bf U} - z\begin{bmatrix}
0 & {\bf 0} \\
{\bf 0} & {\bf I_{n-1}}
\end{bmatrix}
\right)
\]
are precisely ${1\over \overline{u_{jj}}}$ for $2\le j\le n$ and thus $$p_{\bf U}(z)=u_{11}
\prod_{j=2}^{n}\left({1\over \overline{u_{jj}}}-z\right).$$
In particular, the eigenvalues of the matrix
\[ \begin{pmatrix}
u_{22} & \cdots & u_{2n}\\
 &\ddots& \vdots \\
{\bf 0} & & u_{nn}
\end{pmatrix} - {1\over u_{11}}\begin{pmatrix}u_{21}\\ u_{31} \\\vdots \\ u_{n1}\end{pmatrix}\cdot \begin{pmatrix}u_{12},u_{13},\ldots,u_{1n} \end{pmatrix} \]
are precisely ${1\over \overline{u_{jj}}}$ for $2\le j\le n$.
\end{coro}
\begin{proof}
Since ${\bf U}\in U(n)$, we have $|u_{11}|=\left|\prod_{j=2}^n u_{jj}\right|$ (cf. \cite{Zh11}).
Thus, the assumption $u_{11}\neq 0$ implies that $u_{jj}\neq 0$ for $2\le j\le n$.
By Theorem \ref{Thm_Existence1} and the constructions, there is a holomorphic isometry $f=(f_1,g_1,\ldots,g_{n-1}):(\Delta,g_\Delta)\to (\Delta,g_\Delta)\times (\mathbb B^{n-1},g_{\mathbb B^{n-1}})$ such that $f_1$ is non-constant and $R(f_1(w))=w$, where $R:\mathbb P^1\to \mathbb P^1$ is the rational function defined by
$R(z) = {z p_{\bf U}(z)\over \prod_{j=2}^n(u_{jj}-z)}$.
By the proof of Lemma \ref{lem:RationalF1}, we have
\[ R(z) = u_{11} z
\prod_{j=2}^{n}{z-{1\over \overline{u_{jj}}}\over z- u_{jj}} 
= u_{11} z
\prod_{j=2}^{n}{{1\over \overline{u_{jj}}}-z\over u_{jj}-z}.
\]
and thus $p_{\bf U}(z)=u_{11}
\prod_{j=2}^{n}\left({1\over \overline{u_{jj}}}-z\right)$.
Note that
\[ p_{\bf U}(z) = u_{11}\cdot \det\left( \begin{pmatrix}
u_{22} & \cdots & u_{2n}\\
 &\ddots& \vdots \\
{\bf 0} & & u_{nn}
\end{pmatrix} - {1\over u_{11}}\begin{pmatrix}u_{21}\\u_{31}\\\vdots \\ u_{n1}\end{pmatrix}\cdot \begin{pmatrix}u_{12},u_{13},\ldots, u_{1n} \end{pmatrix} - z {\bf I_{n-1}} \right) \]
by item $(2)$ of Remark \ref{remark2.4}.
The proof is complete.
\end{proof}
Actually, Lemma \ref{lem:RationalF1} implies the following theorem due to the fact that $0<|u_{jj}|\le 1$ while $\left|{1\over \overline{u_{jj}}}\right|^2\ge 1$ for $2\le j\le n+1$.

\begin{theo}\label{Thm_Reduction1}
Let $f=(f_1,f_{2,1},\ldots,f_{2,n}): (\Delta,g_\Delta)\to (\Delta,g_\Delta)\times(\mathbb B^n,g_{\mathbb B^n})$ be a holomorphic isometry such that $f_1$ is non-constant, where $n\ge 2$ is an integer.
Let $R:\mathbb P^1\to \mathbb P^1$ be the rational function such that $R(f_1(w))=w$.
Then, we have the following:
\begin{enumerate}
\item[(1)] If $\deg(R)=m+1\le n$ for some positive integer $m$, then $f$ is congruent to a map $\widetilde f: (\Delta,g_\Delta)\to (\Delta,g_\Delta)\times(\mathbb B^n,g_{\mathbb B^n})$ defined by $\widetilde f(w)=(F(w),{\bf 0})$ for some holomorphic isometry $F:(\Delta,g_\Delta)\to (\Delta,g_\Delta)\times(\mathbb B^{m},g_{\mathbb B^{m}})$ such that $F$ is not congruent to
$(\widetilde F,{\bf 0})$ for any holomorphic isometry $\widetilde F:(\Delta,g_\Delta)\to (\Delta,g_\Delta)\times(\mathbb B^{m'},g_{\mathbb B^{m'}})$ with $m'\le m-1$ \rm{(resp. $\widetilde F:(\Delta,g_\Delta)\to (\Delta,g_\Delta)$)}
whenever $m\ge 2$ \rm{(resp. $m=1$)}.
\item[(2)]
If $\deg(R)=1$, then $f$ is congruent to a map $\widetilde f: (\Delta,g_\Delta)\to (\Delta,g_\Delta)\times(\mathbb B^n,g_{\mathbb B^n})$ defined by $\widetilde f(w)=(F(w),{\bf 0})$ for some holomorphic isometry $F:(\Delta,g_\Delta)\to (\Delta,g_\Delta)$.
\end{enumerate}
\end{theo}
\begin{proof}
We first prove $(1)$.
Assume without loss of generality that $f(0)={\bf 0}$.
Then, we have shown that there is a unitary matrix ${\bf U}:=\begin{pmatrix}
u_{ij}
\end{pmatrix}_{1\le i,j\le n+1}\in U(n+1)$ such that
\[ {\bf U} \begin{pmatrix}
f_1(w),f_{2,1}(w),\ldots, f_{2,n}(w)
\end{pmatrix}^T
= \begin{pmatrix}
w,f_1(w)  f_{2,1}(w),\ldots, f_1(w) f_{2,n}(w)
\end{pmatrix}^T. \]
We may further assume without loss of generality that $\begin{pmatrix}
u_{ij}
\end{pmatrix}_{2 \le i,j\le n+1}$ is upper triangular (cf. Section \ref{Sec:2.1.1}).
From the construction and Lemma \ref{lem:RationalF1}, the rational function $R$ is given by
\[ R(z) = u_{11} z
\prod_{j=2}^{n+1}{z-{1\over \overline{u_{jj}}}\over z- u_{jj}}. \]
If $\deg(R)=m+1\le n$, then we have ${1\over \overline{u_{j_\mu j_\mu}}}=u_{l_\mu l_\mu}$ for some $j_\mu,l_\mu$, $2\le j_\mu,l_\mu\le n+1$, and $1\le \mu \le n-m$ such that $l_1,\ldots,l_{n-m}$ are distinct.
But then this implies that for each $\mu$, $1\le \mu\le n-m$, $|u_{l_\mu l_\mu}|^2=1$ so that
$u_{l_\mu k}=u_{k l_\mu}=0$ for any $k\neq l_\mu$.
This would force $f_{2,l_\mu-1}\equiv 0$ for $1\le \mu\le n-m$.
This shows that $f$ is congruent to a map $\widetilde f: (\Delta,g_\Delta)\to (\Delta,g_\Delta)\times(\mathbb B^n,g_{\mathbb B^n})$ defined by $\widetilde f(w)=(F(w),{\bf 0})$ for some holomorphic isometry $F:(\Delta,g_\Delta)\to (\Delta,g_\Delta)\times(\mathbb B^{m},g_{\mathbb B^{m}})$.
Since $\deg(R)=m+1$, we see that it is impossible that $F$ is congruent to
$(\widetilde F,{\bf 0})$ for some holomorphic isometry $\widetilde F:(\Delta,g_\Delta)\to (\Delta,g_\Delta)\times(\mathbb B^{m'},g_{\mathbb B^{m'}})$ for $m'\le m-1$ (resp. $\widetilde F:(\Delta,g_\Delta)\to (\Delta,g_\Delta)$)
whenever $m\ge 2$ (resp. $m=1$).
The case where $\deg(R)=1$ is also clear by our arguments and thus the assertion of $(2)$ follows.
\end{proof}
\begin{rema}
In the settings of Theorem \ref{Thm_Reduction1}, if $\deg(R)=2$, then $f$ is congruent to the holomorphic isometry
$\widetilde f:(\Delta,g_\Delta)\to (\Delta,g_\Delta)\times(\mathbb B^{n},g_{\mathbb B^{n}})$ defined by $\widetilde f(w):=(\alpha(w);\beta(w),0,\ldots,0)$, where $(\alpha,\beta)$ $:$ $(\Delta,ds_\Delta^2)$ $\to$ $(\Delta^2,ds_{\Delta^2}^2)$ is the square-root embedding (cf. \cite{Ng10}).
\end{rema}

\begin{lemm}\label{lem_ExistenceUM1}
For any integer $n\ge 3$, there is a unitary matrix ${\bf U}=\begin{pmatrix}
u_{ij}
\end{pmatrix}_{1\le i,j\le n}\in U(n)$ such that
$u_{kj}=0$ for $2\le j\le k-1$, $3\le k\le n$ and $0<|u_{ll}|<1$ for $2\le l \le n$.
\end{lemm}
\begin{proof}
For $n=3$, we have constructed an explicit matrix in $U(3)$ which satisfies the desired property, namely the matrix
\[ \begin{pmatrix}
-{1\over 2} & {1\over \sqrt{2}} & {1\over 2}\\
{1\over 2} & {1\over \sqrt{2}} & -{1\over 2}\\
{1\over \sqrt{2}} & 0 & {1\over \sqrt{2}} 
\end{pmatrix} \in U(3). \]
Assume that the statement is true for some integer $m\ge 3$, i.e., there is a unitary matrix ${\bf V}=\begin{pmatrix}
v_{ij}
\end{pmatrix}_{1\le i,j\le m}\in U(m)$ such that
$v_{kj}=0$ for $2\le j\le k-1$, $3\le k\le m$ and $0<|v_{jj}|<1$ for $2\le j \le m$.
Then, we see that the matrix
\[ {\bf V}':=\begin{pmatrix}
v_{21} & 0 & v_{22}  & \cdots & v_{2m}\\
\vdots & \vdots & \vdots  & \ddots & \vdots \\
v_{m1} & 0 & 0 & \cdots & v_{mm}
\end{pmatrix}\in M(m-1,m+1;\mathbb C) \]
satisfies ${\bf V'}\overline{\bf V'}^T = {\bf I_{m-1}}$.

Let $W$ be the $\mathbb C$-linear span of the row vectors $\begin{pmatrix}
0 ,1,0,\ldots,0
\end{pmatrix}$, $\begin{pmatrix}
v_{11},0 ,v_{12},\ldots,v_{1m}
\end{pmatrix}$ $\in$ $M(1,m+1;\mathbb C)$.
Then, we have
$\widetilde{\bf U}:=\begin{bmatrix}
{\bf u}_1^T& {\bf u}_2^T & {\bf V'}^T
\end{bmatrix}^T \in U(m+1)$ for any orthonormal basis $\{{\bf u}_1,{\bf u}_2\}$ of $W$.
Let ${\bf u} = \begin{pmatrix}
{1\over \sqrt{2}} v_{11},{1\over \sqrt{2}},{1\over \sqrt{2}} v_{12}, \ldots, {1\over \sqrt{2}} v_{1m}
\end{pmatrix}$ be a vector in $W$. Then, $\lVert {\bf u} \rVert^2= {1\over 2} \sum_{j=1}^m |v_{1j}|^2 + {1\over 2} =1$.
Thus, we can let ${\bf u}_2={\bf u}$ and ${\bf u}_1$ be a unit vector in $W$ such that ${\bf u}_1$ and ${\bf u}_2$ are orthogonal.
Actually, one may choose 
\[ {\bf u_1}=\begin{pmatrix}
-{1\over \sqrt{2}} v_{11},{1\over \sqrt{2}},-{1\over \sqrt{2}} v_{12}, \ldots, -{1\over \sqrt{2}} v_{1m}
\end{pmatrix}, \]
but this does not affect our arguments.
Then, $\widetilde{\bf U}\in U(m+1)$ satisfies the desired property.
By induction, the proof is complete.
\end{proof}
\noindent By Lemma \ref{lem_ExistenceUM1} and Theorem \ref{Thm_Existence1}, we have the following:
\begin{prop}\label{Pro:Exist_top_degree}
For any integer $n\ge 2$, there is a holomorphic isometry $f=(f_1,f_{2,1},\ldots,f_{2,n}): (\Delta,g_\Delta)\to (\Delta,g_\Delta)\times(\mathbb B^n,g_{\mathbb B^n})$ such that $f_1$ is a non-constant function and $R(f_1(w))=w$ on the unit disk $\Delta$ for some rational function $R:\mathbb P^1\to \mathbb P^1$ of degree $n+1$.
\end{prop}
\begin{proof}
By Lemma \ref{lem_ExistenceUM1}, there is a unitary matrix ${\bf U}=\begin{pmatrix}
u_{ij}
\end{pmatrix}_{1\le i,j\le n+1}\in U(n+1)$ such that
$u_{kj}=0$ for $2\le j\le k-1$, $3\le k\le n+1$ and $0<|u_{ll}|<1$ for $2\le l \le n+1$.
Then, it follows from Theorem \ref{Thm_Existence1} that there is a holomorphic isometry $f=(f_1,f_{2,1},\ldots,f_{2,n}): (\Delta,g_\Delta)\to (\Delta,g_\Delta)\times(\mathbb B^n,g_{\mathbb B^n})$ such that $f_1$ is a non-constant function,
\[ {\bf U} \begin{pmatrix}
f_1(w),f_{2,1}(w),\ldots, f_{2,n}(w)
\end{pmatrix}^T
= \begin{pmatrix}
w, f_1(w)  f_{2,1}(w),\ldots,f_1(w) f_{2,n}(w)
\end{pmatrix}^T \]
and there is a rational function $R:\mathbb P^1\to \mathbb P^1$ satisfying $R(f_1(w))=w$.
From the construction, we see that
\[ R(z) =u_{11} z \prod_{j=2}^{n+1} {z-{1\over \overline{u_{jj}}}\over z-u_{jj}}. \]
Since $0<|u_{ll}|<1$ for $2\le l \le n+1$, ${1\over \overline{u_{jj}}}\neq u_{ll}$ for $2\le l,j\le n+1$.
Thus, $\deg(R)=n+1$ and we are done.
\end{proof}
It is natural to ask whether one can relate any holomorphic isometry $(\Delta,g_\Delta)\to (\Delta,g_\Delta)\times(\mathbb B^n,g_{\mathbb B^n})$, $n\ge 2$, to some holomorphic isometry $(\Delta,g_\Delta)\to (\Delta,g_\Delta)\times(\mathbb B^{n-1},g_{\mathbb B^{n-1}})$ with extra parameters.
The following yields certain relations between any given holomorphic isometry $f:(\Delta,g_\Delta)\to (\Delta,g_\Delta)\times(\mathbb B^n,g_{\mathbb B^n})$, $n\ge 2$, and some holomorphic isometry $\widetilde f:(\Delta,g_\Delta)\to (\Delta,g_\Delta)\times(\mathbb B^{n-1},g_{\mathbb B^{n-1}})$ together with one extra parameter.
\begin{theo}\label{Thm:Red_Para1}
Let $f=(f_1,f_{2,1},\ldots,f_{2,n}): (\Delta,g_\Delta)\to (\Delta,g_\Delta)\times(\mathbb B^n,g_{\mathbb B^n})$ be a holomorphic isometry such that $f_1$ is non-constant, where $n\ge 2$ is an integer.
Then, $f$ can be determined by some holomorphic isometry $\widetilde f:(\Delta,g_\Delta)\to (\Delta,g_\Delta)\times(\mathbb B^{n-1},g_{\mathbb B^{n-1}})$ and some parameter $\zeta\in \overline\Delta\smallsetminus\{0\}$ up to congruence.
\end{theo}
\begin{proof}
Assume without loss of generality that $f(0)={\bf 0}$.
By \cite[Theorem 2]{Ca53} we have
\[  {\bf U} \begin{pmatrix}
f_1(w),f_{2,1}(w),\ldots, f_{2,n}(w)
\end{pmatrix}^T
= \begin{pmatrix}
w, f_1(w)  f_{2,1}(w),\ldots,f_1(w) f_{2,n}(w)
\end{pmatrix}^T \]
for some ${\bf U}\in U(n+1)$ and for any $w\in \Delta$.
Then, we have obtained another holomorphic isometry $F$ $:=$ $(f_1,g_1,\ldots,g_n):(\Delta,g_\Delta)\to (\Delta,g_\Delta)\times(\mathbb B^n,g_{\mathbb B^n})$ such that $F$ is congruent to $f$, $F(0)={\bf 0}$ and
\[{\bf U}' \begin{pmatrix}
f_1(w),g_{1}(w) ,\ldots, g_n(w)
\end{pmatrix}^T
= \begin{pmatrix}
w, f_1(w)  g_{1}(w) ,\ldots, f_1(w) g_{n}(w)
\end{pmatrix}^T\,\,\forall\;w\in \Delta,\]
where ${\bf U}'=\begin{pmatrix}
u_{ij}
\end{pmatrix}_{1\le i,j\le n+1}\in U(n+1)$ is some unitary matrix such that
$\begin{pmatrix}
u_{ij}
\end{pmatrix}_{2\le i,j\le n+1}$ is an upper triangular matrix.
This does not affect the rational function $R:\mathbb P^1 \to \mathbb P^1$ which satisfies $R(f_1(w))=w$. If $\deg(R)\le n$, then we are done by Theorem \ref{Thm_Reduction1}.
Therefore, we now consider the case where $\deg(R)=n+1$ so that $0<|u_{jj}|<1$ for $2\le j\le n+1$.
Then, we obtain the matrix
\[ \widetilde{\bf U}:= \begin{pmatrix}
a_1 & a_2 & \cdots & a_n\\
u_{31} & u_{33} &\cdots & u_{3,n+1}\\
\vdots & \vdots & \ddots & \vdots\\
u_{n+1,1} & 0 & \cdots & u_{n+1,n+1}
\end{pmatrix}\in U(n) \]
for some row vector $\begin{pmatrix}
a_1,\ldots,a_n
\end{pmatrix}\in M(1,n;\mathbb C)$.
Note that we have $|a_1|= \prod_{j=3}^{n+1} |u_{jj}|$.
Then, one observes that the row vectors
${\bf u}_j:=\begin{pmatrix}
u_{j1},\ldots,u_{j,n+1}
\end{pmatrix}$, $j=1,2$, are $\mathbb C$-linear combinations of $\begin{pmatrix}
a_1,0,a_2,\ldots,a_n
\end{pmatrix}$ and $\begin{pmatrix}
0,1,0,\ldots,0
\end{pmatrix}$.
In particular, we have
${\bf u}_2= \begin{pmatrix}c_1 a_1,c_2,c_1a_2,\ldots,c_1a_n \end{pmatrix}$
for some $c_1,c_2\in \mathbb C$ such that $|c_1|^2+|c_2|^2=1$.
Now, we have $u_{22}=c_2$ so that $0<|c_2|\le 1$ from the construction.
Moreover, this determines that
\[ {\bf u}_1 = e^{i\theta} \begin{pmatrix}\overline{c_2} a_1, -\overline{c_1}, \overline{c_2} a_2,\ldots,\overline{c_2} a_n\end{pmatrix} \]
for some $\theta\in [0,2\pi)$.
By Theorem \ref{Thm_Existence1}, the unitary matrix $\widetilde{\bf U}$ defines a holomorphic isometry $\widetilde f=(\widetilde f_1,\widetilde g_1,\ldots,\widetilde g_{n-1}):(\Delta,g_\Delta)\to (\Delta,g_\Delta)\times(\mathbb B^{n-1},g_{\mathbb B^{n-1}})$ and a rational function $\widetilde R:\mathbb P^1\to \mathbb P^1$ such that $\widetilde R(\widetilde f_1 (w))=w$.
From the constructions, we have
\begin{equation}\label{Eq:RatDecom}
R(z) = e^{i\theta} \overline{c_2} \widetilde R(z) {z-{1\over \overline{c_2}}\over z- c_2} 
= e^{i\theta}  \widetilde R(z) {\overline{c_2}z-1\over z- c_2} 
\end{equation}
for some $\theta\in [0,2\pi)$.
(Noting that $\deg(\widetilde R) = n$ and $0<|c_2|<1$ under the assumption that $\deg(R)=n+1$.)
We may assume without loss of generality that $\theta=0$ because this does not affect the equivalence class of $f$.
Actually, Eq. (\ref{Eq:RatDecom}) is still valid when $\deg(R)\le n$.

It follows that the rational function $R$ is determined by the holomorphic isometry $\widetilde f:(\Delta,g_\Delta)\to (\Delta,g_\Delta)\times(\mathbb B^{n-1},g_{\mathbb B^{n-1}})$ and a parameter $c_2\in \overline\Delta\smallsetminus\{0\}$.
From item $(3)$ of Remark \ref{remark2.4}, $R$ determines $f$ uniquely up to congruence so that $f$ actually depends on the holomorphic isometry $\widetilde f:(\Delta,g_\Delta)\to (\Delta,g_\Delta)\times(\mathbb B^{n-1},g_{\mathbb B^{n-1}})$ and a parameter $c_2\in \overline\Delta\smallsetminus\{0\}$ up to congruence.
\end{proof}
\begin{rema}
\begin{enumerate}
\item[(1)]
From the above theorem, given a holomorphic isometry $\widetilde f:(\Delta,g_\Delta)\to (\Delta,g_\Delta)\times(\mathbb B^{n-1},g_{\mathbb B^{n-1}})$ and a parameter $\zeta\in \overline\Delta\smallsetminus\{0\}$, we have a holomorphic isometry $F_{\widetilde f,\zeta}: (\Delta,g_\Delta)\to (\Delta,g_\Delta)\times(\mathbb B^{n},g_{\mathbb B^{n}})$.
It is not known whether $F_{\widetilde f,\zeta}$ and $F_{\widetilde f,\zeta'}$ are congruent or incongruent to each other for distinct $\zeta,\zeta'\in \Delta\smallsetminus\{0\}$ in general.
\item[(2)] For any $\zeta\in \partial\Delta$ and any holomorphic isometry $\widetilde f=(\widetilde f_1,\widetilde g_1,\ldots,\widetilde g_{n-1}):(\Delta,g_\Delta)\to (\Delta,g_\Delta)\times(\mathbb B^{n-1},g_{\mathbb B^{n-1}})$, the holomorphic isometry $F_{\widetilde f,\zeta}: (\Delta,g_\Delta)\to (\Delta,g_\Delta)\times(\mathbb B^{n},g_{\mathbb B^{n}})$ is congruent to
the map $($$\widetilde f_1$,\,$0$, $\widetilde g_1$, $\ldots$, $\widetilde g_{n-1}$$)$.
Thus, $F_{\widetilde f,\zeta}$ and $F_{\widetilde f,\zeta'}$ are congruent to each other for any distinct $\zeta,\zeta'\in \partial\Delta$.
\end{enumerate}
\end{rema}

\subsection{The case where the target is $\Delta\times \mathbb B^2$}\label{Sec:1.1}
The following is a corollary of Theorem \ref{Thm_Reduction1}.
In addition, we can characterize those holomorphic isometries $(\Delta,g_\Delta)\to (\Delta,g_\Delta)\times(\mathbb B^2,g_{\mathbb B^2})$ obtained from the square-root embedding $(\Delta,ds_\Delta^2)\to(\Delta^2,ds_{\Delta^2}^2)$ (cf. Ng \cite{Ng10}).
\begin{coro}
\label{Cor:Char_Square-root1}
Let $f=(f_1,f_{2,1},f_{2,2}):(\Delta,g_\Delta)\to (\Delta,g_\Delta)\times(\mathbb B^2,g_{\mathbb B^2})$ be a holomorphic isometry.
Suppose that $f_1$ is a non-constant function.
Then, $\deg(R)\le 2$ if and only if $f$ is congruent to either $w\mapsto (w;0,0)$ or $w\mapsto(\alpha_1(w);\beta_1(w),0)$, where $(\alpha_1,\beta_1):(\Delta,ds_\Delta^2)\to(\Delta^2,ds_{\Delta^2}^2)$ is the square-root embedding \rm{(cf. Ng \cite{Ng10})}, and $R:\mathbb P^1\to \mathbb P^1$ is the rational function such that $R(f_1(w))=w$.
\end{coro}
\begin{proof}
If $f$ is congruent to one of the given holomorphic isometries, then it is clear that $\deg(R)\le 2$.
Conversely, if $\deg(R)\le 2$, then it follows from Theorem \ref{Thm_Reduction1} that $f$ is obtained from some holomorphic isometry $(\Delta,g_\Delta)\to (\Delta,g_\Delta)\times (\Delta,g_\Delta)$ and the rest follows from the classification for the $2$-disk by Ng \cite{Ng10}.
\end{proof}
\begin{rema}
\begin{enumerate}
\item[(1)] Actually, we have $\deg(R)=2$ if and only if the holomorphic isometry $f:(\Delta,g_\Delta)\to (\Delta,g_\Delta)\times(\mathbb B^2,g_{\mathbb B^2})$ is congruent to the map $(\Delta,g_\Delta)\to (\Delta,g_\Delta)\times(\mathbb B^2,g_{\mathbb B^2})$ given by $w\mapsto(\alpha_1(w);\beta_1(w),0)$, where $(\alpha_1,\beta_1):(\Delta,ds_\Delta^2)\to(\Delta^2,ds_{\Delta^2}^2)$ is the square-root embedding.
\item[(2)] It follows from Corollary \ref{Cor:Char_Square-root1} and Proposition \ref{Pro:Exist_top_degree} that there exists a holomorphic isometry $f=(f_1,g_1,g_2):(\Delta,g_\Delta)\to (\Delta,g_\Delta)\times(\mathbb B^2,g_{\mathbb B^2})$ such that $\deg(R)=3$ and $R(f_1(w))=w$ for some rational function $R:\mathbb P^1 \to \mathbb P^1$ so that $f$ is incongruent to any of the holomorphic isometris $(\Delta,g_\Delta)\to (\Delta,g_\Delta)\times(\mathbb B^2,g_{\mathbb B^2})$ given by $w\mapsto (w;0,0)$, $w\mapsto(0;w,0)$ or $w\mapsto(\alpha_1(w);\beta_1(w),0)$, where $(\alpha_1,\beta_1):(\Delta,ds_\Delta^2)\to(\Delta^2,ds_{\Delta^2}^2)$ is the square-root embedding.
\end{enumerate}
\end{rema}

\subsection{Existence of a real $1$-parameter family of mutually incongruent holomorphic isometries and generalizations}
\text{}\\
Let $f=(f_1,g_1,\ldots,g_n):(\Delta,g_\Delta)\to (\Delta,g_\Delta)\times (\mathbb B^n,g_{\mathbb B^n})$ and
$\widetilde f=(\widetilde f_1,\widetilde g_1$, $\ldots$, $\widetilde g_n):(\Delta,g_\Delta)\to (\Delta,g_\Delta)\times (\mathbb B^n,g_{\mathbb B^n})$ be holomorphic isometries such that $f_1$ and $\widetilde f_1$ are non-constant functions, where $n\ge 2$ is an integer.
We may suppose that $f(0)=\widetilde f(0)={\bf 0}$ without loss of generality.
Then, there are rational functions $R,\widetilde R:\mathbb P^1\to \mathbb P^1$ such that $R(f_1(w))=w$ and $\widetilde R(\widetilde f_1(w))=w$.
If $f$ is congruent to $\widetilde f$, then we have $\psi\circ f_1 \circ \varphi= \widetilde f_1$ for some $\varphi,\psi\in \mathrm{Aut}(\Delta)$ so that
\[ R = \varphi\circ\widetilde R\circ\psi. \]
In particular, $\psi$ maps the ramification locus of $R$ onto the ramification locus of $\widetilde R$.
In addition, $\varphi$ maps the branch locus of $\widetilde R$ onto the branch locus of $R$.

A slight modification of the proof of Lemma \ref{lem_ExistenceUM1} yields the following:
\begin{lemm}\label{lem_ExistenceUM2}
For any integer $n\ge 3$, there is a unitary matrix ${\bf U}=\begin{pmatrix}
u_{ij}
\end{pmatrix}_{1\le i,j\le n}\in U(n)$ such that
$u_{kj}=0$ for $2\le j\le k-1$, $3\le k\le n$, $u_{22}=\cdots = u_{nn} =:\zeta\in \Delta\smallsetminus\{0\}$.
\end{lemm}
\begin{proof}
For $n=3$, we have the matrix
\[ {\bf U}_{\zeta}:=\begin{pmatrix}
-\overline{\zeta}^2 & -\sqrt{1-|\zeta|^2} & \overline{\zeta} \sqrt{1-|\zeta|^2} \\
- \sqrt{1-|\zeta|^2} \cdot \overline{\zeta} & \zeta & 1-|\zeta|^2 \\
\sqrt{1-|\zeta|^2} & 0 & \zeta
\end{pmatrix}\in U(3) \]
which satisfies the desired properties.
For $m\ge 3$, we simply let 
\[ {\bf u}_2 = \begin{pmatrix}
\sqrt{1-|\zeta|^2} v_{11},\zeta, \sqrt{1-|\zeta|^2}  v_{12}, \ldots, \sqrt{1-|\zeta|^2}  v_{1m}
\end{pmatrix}\]
(resp. ${\bf u}_1=\begin{pmatrix}
-\overline{\zeta} v_{11},\sqrt{1-|\zeta|^2} ,-\overline{\zeta} v_{12}, \ldots, -\overline{\zeta} v_{1m}
\end{pmatrix}$)
in place of the original ${\bf u}_2$ (resp. ${\bf u}_1$) in the proof of Lemma \ref{lem_ExistenceUM1}.
Then, we also obtain $\widetilde{\bf U}\in U(m+1)$ as in the proof of Lemma \ref{lem_ExistenceUM1} and $\widetilde{\bf U}$ satisfies the desired properties.
The proof is complete by induction.
\end{proof}
The following shows the existence of a real $1$-parameter family $\{f_t\}_{t\in \mathbb R}$ of mutually incongruent holomorphic isometries $f_{t}:(\Delta,g_\Delta)\to (\Delta,g_\Delta)\times (\mathbb B^n,g_{\mathbb B^n})$.
\begin{prop}\label{Pro:Ex_1PFICHIndim}
Let $n\ge 2$ be an integer.
Then, there is a real $1$-parameter family $\{f_t\}_{t\in \mathbb R}$ of mutually incongruent holomorphic isometries $f_{t}:(\Delta,g_\Delta)\to (\Delta,g_\Delta)\times (\mathbb B^n,g_{\mathbb B^n})$.

More generally, there is a family $\{f_\zeta\}_{\zeta\in A_n}$ of holomorphic isometries $f_{\zeta}:(\Delta,g_\Delta)\to (\Delta,g_\Delta)\times (\mathbb B^n,g_{\mathbb B^n})$ such that for any $\zeta,\zeta'\in A_n:=\left\{\xi\in \mathbb C\mid {n-1\over n+1}<|\xi|<1\right\}$, $f_\zeta$ and $f_{\zeta'}$ are congruent to each other if and only if $|\zeta|=|\zeta'|$.

In addition, if $n\ge 3$, then there is a family $\{f_\zeta\}_{\zeta\in \Delta^*}$ of holomorphic isometries $f_{\zeta}:(\Delta,g_\Delta)\to (\Delta,g_\Delta)\times (\mathbb B^n,g_{\mathbb B^n})$ such that for any $\zeta,\zeta'\in \Delta^*:=\Delta\smallsetminus\{0\}$, $f_\zeta$ and $f_{\zeta'}$ are congruent to each other if and only if $|\zeta|=|\zeta'|$.
\end{prop}
\begin{proof}
For any integer $n\ge 2$, there is a matrix ${\bf U}_\zeta = \begin{pmatrix}
u_{ij}(\zeta)
\end{pmatrix}_{1\le i,j\le n+1}$ $\in$ $U(n+1)$ such that
$u_{kj}(\zeta)=0$ for $2\le j\le k-1$, $3\le k\le n+1$, $u_{jj}(\zeta)=\zeta$ for $2\le j \le n+1$ and $u_{11}(\zeta)=\overline\zeta^n$ for any $\zeta\in \Delta\smallsetminus\{0\}$ by Lemma \ref{lem_ExistenceUM2}.

By Theorem \ref{Thm_Existence1}, there is a holomorphic isometry $f_\zeta=(f_\zeta^1,f_\zeta^{2,1},\ldots,f_\zeta^{2,n}):(\Delta,g_\Delta)\to (\Delta,g_\Delta)\times (\mathbb B^n,g_{\mathbb B^n})$ such that $f_\zeta(0)={\bf 0}$ and
\[ {\bf U}_{\zeta} \begin{pmatrix}
f_\zeta^1(w),
f_\zeta^{2,1}(w),\ldots,
f_\zeta^{2,n}(w)
\end{pmatrix}^T
=\begin{pmatrix}
w,
f_\zeta^1(w)f_\zeta^{2,1}(w),\ldots,f_\zeta^1(w)f_\zeta^{2,n}(w)
\end{pmatrix}^T \]
for any $w\in \Delta$.
Moreover, we have $R_\zeta(f_\zeta^1(w))=w$, where $R_\zeta:\mathbb P^1\to \mathbb P^1$ is the rational function of degree $n+1$ given by
\[ R_\zeta(z)=z \left({\overline\zeta z-1\over z-\zeta}\right)^n. \]
We observe that
\[ \begin{split}
R_{\zeta e^{i\theta}}(z)
=& z \left({\overline\zeta e^{-i\theta}z-1\over z-e^{i\theta}\zeta}\right)^n
= e^{-i(n-1)\theta}(e^{-i\theta}z)\left({\overline\zeta e^{-i\theta}z-1\over e^{-i\theta}z-\zeta}\right)^n\\
=& e^{-i(n-1)\theta} R_\zeta(e^{-i\theta}z) 
\end{split}\]
for any $\theta\in [0,2\pi)$.
Since $R_\zeta(e^{-i\theta} f_{\zeta e^{i\theta}}^1(e^{-i(n-1)\theta}w))=w$ and $0$ is not a ramification point (or critical point) of $R_\zeta$, we have
$$f_\zeta^1(w) = e^{-i\theta} f_{\zeta e^{i\theta}}^1(e^{-i(n-1)\theta}w)$$ for any $w\in \Delta$
by the Identity Theorem of holomorphic functions.
In particular, $1-|f_\zeta^1(w)|^2 = 1-| f_{\zeta e^{i\theta}}^1(e^{-i(n-1)\theta}w)|^2$ for any $w\in \Delta$.
Then, it follows from the functional equations
\[ \begin{split}
&(1-|f_\zeta^1(w)|^2)\left(1-\sum_{j=1}^n |f_\zeta^{2,j}(w)|^2\right)\\
=&1-|w|^2=1-|e^{-i(n-1)\theta}w|^2\\
=& (1-| f_{\zeta e^{i\theta}}^1(e^{-i(n-1)\theta}w)|^2)
\left(1-\sum_{j=1}^n |f_{\zeta e^{i\theta}}^{2,j}(e^{-i(n-1)\theta}w)|^2\right) \end{split}\]
for any $w\in \Delta$
that
$\sum_{j=1}^n |f_\zeta^{2,j}(w)|^2 = \sum_{j=1}^n |f_{\zeta e^{i\theta}}^{2,j}(e^{-i(n-1)\theta} w)|^2$ holds on $\Delta$.
Then, there exists a unitary transformation $\Psi:\mathbb C^n \to \mathbb C^n$ such that
$$\Psi(f_\zeta^{2}(w)) = f_{\zeta e^{i\theta}}^{2}(e^{-i(n-1)\theta} w)$$ for any $w\in \Delta$ by the local rigidity theorem of Calabi \cite[Theorem 2]{Ca53}, where
$f_\chi^2:=(f_\chi^{2,1},\ldots,f_\chi^{2,n})$ for any $\chi\in \Delta\smallsetminus\{0\}$.
Actually, $\Psi$ belongs to the isotropy subgroup of $\mathrm{Aut}(\mathbb B^n)$ at ${\bf 0}$, and the map $\phi$ defined by $\phi(w):= e^{-i(n-1)\theta} w$ belongs to the isotropy subgroup of $\mathrm{Aut}(\Delta)$ at $0$.
Let $\Phi:\mathbb C\times \mathbb C^n\to \mathbb C\times \mathbb C^n$ be the map
\[ \Phi(z,w_1,\ldots,w_n)=(e^{i\theta} z, \Psi(w_1,\ldots,w_n)). \]
Then, it is clear that $\Phi$ belongs to the isotropy subgroup of the automorphism group of $\Delta\times \mathbb B^n$ at $(0,{\bf 0})$ and
\[ \Phi(f_\zeta(w))=f_{\zeta e^{i\theta}}(e^{-i(n-1)\theta}w) 
=f_{\zeta e^{i\theta}}(\phi(w)) \quad \forall \;w\in \Delta\]
so that
$\Phi\circ f_\zeta \circ \phi^{-1} = f_{\zeta e^{i\theta}}$, i.e., $f_\zeta$ and $f_{\zeta e^{i\theta}}$ are congruent to each other.
In particular, if $\zeta,\zeta'\in \Delta^*=\Delta\smallsetminus\{0\}$ such that $|\zeta|=|\zeta'|$, then $f_{\zeta}$ and $f_{\zeta'}$ are congruent to each other.

Now, we will show that $f_{\zeta}$ and $f_{\zeta'}$ are incongruent to each other whenever $|\zeta|\neq |\zeta'|$ for $\zeta,\zeta'\in \Delta\smallsetminus\{0\}$.
Write $R_\zeta(z)={p_\zeta(z)\over q_\zeta(z)}$, where
$p_\zeta(z):=z(\overline \zeta z-1)^n$ and $q_\zeta(z):=(z-\zeta)^n$.
Then, the ramification points of $R_\zeta$ are the zeros of $p_\zeta'(z)q_\zeta(z)-p_\zeta(z)q_\zeta'(z)$.
In particular, the ramification points of $R_\zeta$ are precisely
${1\over \overline\zeta},\zeta,a_+(\zeta),a_-(\zeta)$, where
\[ \begin{split} a_\pm(\zeta)&:=
{(n+1)|\zeta|^2+(1-n) \pm \sqrt{(n-1)^2-(2n^2+2)|\zeta|^2 + (n+1)^2 |\zeta|^4}\over 2\overline\zeta}\\
&\in \mathbb C\smallsetminus\{0\}.\end{split} \]
In addition, the ramification order of $R_\zeta$ at ${1\over \overline\zeta}$ (resp. $\zeta$) is equal to $n-1$ and the ramification order of $R_\zeta$ at $a_+(\zeta)$ (resp. $a_-(\zeta)$) is equal to $1$.
Therefore, $R_\zeta(a_{\pm}(\zeta))\in \mathbb C \smallsetminus \{0\}$.
Note that $a_+(\zeta)=a_-(\zeta)={(n+1)|\zeta|^2+(1-n) \over 2\overline\zeta}$ whenever $\zeta$ satisfies $(n-1)^2-(2n^2+2)|\zeta|^2 + (n+1)^2 |\zeta|^4=0$.
It is obvious that
$(n-1)^2-(2n^2+2)|\zeta|^2 + (n+1)^2 |\zeta|^4<0$ whenever ${n-1\over n+1}<|\zeta|<1$.
Then, we have
\[ a_{\pm}(\zeta)
={(n+1)|\zeta|^2+(1-n) \pm \sqrt{(2n^2+2)|\zeta|^2 - (n+1)^2 |\zeta|^4-(n-1)^2}\;i\over 2\overline\zeta} \]
with $\sqrt{(2n^2+2)|\zeta|^2 - (n+1)^2 |\zeta|^4-(n-1)^2}\in \mathbb R\smallsetminus\{0\}$ for ${n-1\over n+1}<|\zeta|<1$, where $i=\sqrt{-1}$.
Moreover, we have $|a_{\pm}(\zeta)|^2=1$ whenever ${n-1\over n+1}<|\zeta|<1$.

Let $A_n:=\left\{\xi\in \mathbb C: {n-1\over n+1}<|\xi|<1\right\}$, where $n\ge 2$.
Now, we simply write $a_{\pm}(\zeta)=e^{i\theta^\pm_{\zeta}}$ for $\zeta\in A_n$.
Note that $R_\zeta(\partial\Delta)\subset \partial\Delta$ (cf. Lemma \ref{lem:RationalF1}).
For any $\zeta\in A_n$, the branch points of $R_\zeta$ are $0,\infty$,
$R_\zeta(e^{i\theta^+_{\zeta}})=:e^{i\phi^+_\zeta}$ and
$R_\zeta(e^{i\theta^-_{\zeta}})=:e^{i\phi^-_\zeta}$.
Note that a priori it is possible that $R_\zeta(e^{i\theta^+_{\zeta}})=R_\zeta(e^{i\theta^-_{\zeta}})=:e^{i\phi^-_\zeta}$ when $n\ge 3$ because $R_\zeta$ is $(n+1)$-sheeted and the ramification order of $R_\zeta$ at $e^{i\theta^\pm_{\zeta}}$ is equal to $1$.

Given any distinct $\zeta,\zeta'\in A_n$, we suppose that $R_\zeta = \varphi \circ R_{\zeta'} \circ \psi$ for some $\psi,\varphi\in \mathrm{Aut}(\Delta)$.
Then, $\varphi$ maps the branch locus of $R_{\zeta'}$ onto the branch locus of $R_\zeta$.
From the fact that $\varphi\in \mathrm{Aut}(\Delta)$ and the above observations, we have $\varphi(0)=0$ so that $\varphi(w)=e^{i\theta_1}w$ for some $\theta_1\in [0,2\pi)$.
Now, zeros of $R_\zeta(z)$ and $R_{\zeta'}(\psi(z))$ are the same by the assumption.
Note that the zeros of $R_{\xi}$ are $0$ and ${1\over \overline\xi}\in \mathbb C\smallsetminus \overline\Delta$.
Thus, we have $\psi^{-1}(0)=0$ and $\psi^{-1}\left({1\over \overline{\zeta'}}\right) = {1\over \overline\zeta}$ since $\psi\in \mathrm{Aut}(\Delta)$.
In particular, $\psi(w)=e^{i\theta_2}w$ for some $\theta_2\in [0,2\pi)$.
Moreover, $\psi(w)=e^{i\theta_2}w$ and $\psi^{-1}\left({1\over \overline{\zeta'}}\right) = {1\over \overline\zeta}$ implies that $|\zeta'|=|\zeta|$.

If $f_\zeta$ and $f_{\zeta'}$ are congruent to each other, then it is obvious that $R_\zeta = \varphi \circ R_{\zeta'} \circ \psi$ for some $\psi,\varphi\in \mathrm{Aut}(\Delta)$ and thus $|\zeta'|=|\zeta|$ by the above arguments.
In other words, $f_\zeta$ and $f_{\zeta'}$ are incongruent to each other for any $\zeta,\zeta'\in A_n$ such that $|\zeta'|\neq |\zeta|$.
Thus, we have a real $1$-parameter family $\{f_t\}_{t\in \left({n-1\over n+1},1\right)}$ of mutually incongruent holomorphic isometries $f_t:(\Delta,g_\Delta)\to (\Delta,g_\Delta)\times (\mathbb B^n,g_{\mathbb B^n})$.
This finishes the proof of the first statement.
Indeed, the second statement is also proved.

It remains to prove the third statement by showing that $\{f_\zeta\}_{\zeta\in \Delta^*}$ is the desired family of holomorphic isometries when $n\ge 3$.
More precisely, we only need to focus on the case where $0<|\zeta|<{n-1\over n+1}$.
We do not restrict to the case where $n\ge 3$ yet.
Let $B_n:=\left\{\xi\in \mathbb C: 0<|\xi|<{n-1\over n+1}\right\}$ for $n\ge 2$.
Then, we have
\[ a_\pm(\zeta)=
{(n+1)|\zeta|^2+(1-n) \pm \sqrt{(n-1)^2-(2n^2+2)|\zeta|^2 + (n+1)^2 |\zeta|^4}\over 2\overline\zeta} \]
and $\sqrt{(n-1)^2-(2n^2+2)|\zeta|^2 + (n+1)^2 |\zeta|^4}\in \mathbb R$.
In particular,
\[ |a_\pm(\zeta)|^2
= {2(n+1)^2|\zeta|^4 -4n^2|\zeta|^2 + 2(n-1)^2
\pm 
2((n+1)|\zeta|^2+(1-n)) \sqrt{C_\zeta}\over 4|\zeta|^2}, \]
where $C_\zeta:=(n+1)^2|\zeta|^4-2(n^2+1)|\zeta|^2 + (1-n)^2$.
Note that
\[ \begin{split}
&2(n+1)^2|\zeta|^4 -4n^2|\zeta|^2 + 2(n-1)^2
-
2((n+1)|\zeta|^2+(1-n)) \sqrt{C_\zeta} - 4|\zeta|^2 \\
=&
2C_\zeta  - 2((n+1)|\zeta|^2+(1-n)) \sqrt{C_\zeta}
>0
\end{split}\]
because $(n+1)|\zeta|^2+(1-n) < (n+1){(n-1)^2\over (n+1)^2}+1-n
= -{2(n-1)\over n+1}<0$ and $C_\zeta>0$.
Thus, we have $|a_-(\zeta)|^2>1$.
On the other hand, we compute
\[ \begin{split}
&((n+1)|\zeta|^2+(1-n))^2
-\left(\sqrt{(n-1)^2-(2n^2+2)|\zeta|^2 + (n+1)^2 |\zeta|^4}\right)^2\\
=& -2(n^2-1)|\zeta|^2 + (2n^2+2)|\zeta|^2
=4|\zeta|^2,
\end{split}\]
so that 
\[ \begin{split}
&a_+(\zeta)\overline{a_-(\zeta)}\\
=&{((n+1)|\zeta|^2+(1-n))^2-\left(\sqrt{(n-1)^2-(2n^2+2)|\zeta|^2 + (n+1)^2 |\zeta|^4}\right)^2\over 4|\zeta|^2}\\
=&1, 
\end{split}\]
i.e., $a_-(\zeta)={1\over \overline{a_+(\zeta)}}$.
Thus, $|a_+(\zeta)|= {1\over |a_-(\zeta)|}<1$.
In particular, the ramification points of $R_\zeta$ are
${1\over \overline\zeta},\zeta$, $a_+(\zeta)$, $a_-(\zeta)$ with
$0<|a_+(\zeta)|<1$ and $|a_-(\zeta)|>1$ whenever $0<|\zeta|<{n-1\over n+1}$. In addition, none of the ramification points of $R_\zeta$ are on the unit circle $\partial\Delta$.
This implies that $f_\zeta$ and $f_{\zeta'}$ are incongruent to each other for any $\zeta\in B_n$ and $\zeta'\in A_n$.
Note that the branch points of $R_\zeta$ are $0,\infty$, $w_0(\zeta)$ and ${1\over \overline{w_0(\zeta)}}$ for some $w_0(\zeta)\in \overline{\Delta}\smallsetminus\{0\}$.
Here $\left\{w_0(\zeta),{1\over \overline{w_0(\zeta)}}\right\}
=\left\{R_\zeta(a_+(\zeta)),R_\zeta(a_-(\zeta))={1\over \overline{R_\zeta(a_+(\zeta))}}\right\}\subset \mathbb C\smallsetminus\{0\}$.
Note that a priori it is possible that $|w_0(\zeta)|^2=1$ for $n\ge 3$ because $R_\zeta:\mathbb P^1 \to \mathbb P^1$ is $(n+1)$-sheeted and the ramification order of $R_\zeta$ at $a_+(\zeta)$ (resp. $a_-(\zeta)$) is equal to $1$ for $|\zeta|\neq {n-1\over n+1}$.

Given any distinct $\zeta,\zeta'\in B_n$, we suppose that $f_\zeta$ and $f_{\zeta'}$ are congruent to each other.
Then, $f_{\zeta'}^1=\psi\circ f_\zeta^1\circ \varphi$ for some $\varphi,\psi\in \mathrm{Aut}(\Delta)$.
It follows that $R_\zeta = \varphi\circ R_{\zeta'}\circ \psi$.
Then, $\varphi$ maps the branch locus of $R_{\zeta'}$ onto the branch locus of $R_\zeta$.
Moreover, $\varphi$ should preserve the branching order of the branch points.
We will make use of the fact that $\varphi(\Delta)\subset \Delta$ and $\varphi(\mathbb P^1 \smallsetminus \overline\Delta)\subset \mathbb P^1 \smallsetminus \overline\Delta$ as we regard $\varphi$ as an automorphism of $\mathbb P^1$.

Now, we consider the case where $n\ge 3$.
For $\zeta\in B_n$, the branching order of $R_\zeta$ at $0$ is equal to $n-1\ge 2$ and the fiber $R_\zeta^{-1}(0)$ contains a ramification point of ramification order $n-1\ge 2$ while the branching order of $R_\zeta$ at $w_0(\zeta)$ is at most $2$ and the fiber $R_\zeta^{-1}(w_0(\zeta))$ contains ramification point(s) of ramification order $1$.
Thus, it is impossible that $\varphi(w_0(\zeta))=0$ so that $\varphi(0)=0$.
Since $R_{\zeta'}(\psi(0))=0$ and $\psi(0)$ is not a ramification point of $R_{\zeta'}$, we have $\psi(0)=0$ and thus both $\varphi$ and $\psi$ are rotations.
In particular, $\psi(z)=e^{i\theta_2}z$ for some $\theta_2\in [0,2\pi)$ so that ${1\over \overline{\zeta'}}=e^{i\theta_2}{1\over \overline\zeta}$ and $|\zeta'|=|\zeta|$.

Note that $R_\zeta$ only has three distinct ramification points $\zeta,{1\over \overline\zeta}$ and $a_+(\zeta)=a_-(\zeta)$ whenever $|\zeta|={n-1\over n+1}$. Thus $f_{{n-1\over n+1} e^{i\theta}}$ and $f_{\zeta'}$ are incongruent to each other for any $\zeta'\in \Delta\smallsetminus\{0\}$ and $\theta\in [0,2\pi)$ such that $|\zeta'|\neq {n-1\over n+1}$.
Hence, we conclude that for any $\zeta,\zeta'\in \Delta\smallsetminus\{0\}$, $f_\zeta$ and $f_{\zeta'}$ are congruent to each other if and only if $|\zeta|=|\zeta'|$.
\end{proof}
\begin{rema}
Write $\Omega:=\Delta^{n+1} \times \mathbb B^n \times \mathbb B^n$.
From the proof of Proposition \ref{Pro_ExThm} and Proposition \ref{Pro:Ex_1PFICHIndim}, there is a real $1$-parameter family $\{f_t\}_{t\in \mathbb R}$ of mutually incongruent holomorphic isometries $\widetilde f_{t}:(\Delta,(n+1)ds_\Delta^2)\to (\Omega,ds_{\Omega}^2)
$.
In addition, all statements in Proposition \ref{Pro:Ex_1PFICHIndim} are true if we replace $(\Delta,g_\Delta)$ and $(\Delta,g_\Delta)\times (\mathbb B^n,g_{\mathbb B^n})$ by $(\Delta,(n+1)ds_\Delta^2)$ and $(\Omega,ds_{\Omega}^2)$ respectively.
\end{rema}
\begin{coro}\label{Cor:ExHolo_n=2}
In the settings of the proof of Proposition \ref{Pro:Ex_1PFICHIndim}, the holomorphic isometry $f_\zeta=(f_\zeta^1,f_{\zeta}^{2,1},f_{\zeta}^{2,2}):(\Delta,g_\Delta)\to (\Delta,g_\Delta)\times (\mathbb B^2,g_{\mathbb B^2})$ is not totally geodesic, irrational (i.e., some component functions of $f_\zeta=f_\zeta(w)$ are not rational functions of $w\in \Delta\subset \mathbb C$) and extends holomorphically to a neighborhood of the closed unit disk $\overline\Delta$ for $0<|\zeta|<{1\over 3}$.
\end{coro}
\begin{proof}
We fix any $\zeta\in \mathbb C$ such that $0<|\zeta|<{1\over 3}$.
Let $f_\zeta=(f_\zeta^1,f_{\zeta}^{2,1},f_{\zeta}^{2,2}):(\Delta,g_\Delta)\to (\Delta,g_\Delta)\times (\mathbb B^2,g_{\mathbb B^2})$ be the holomorphic isometry constructed in the proof of Proposition \ref{Pro:Ex_1PFICHIndim}.
From the construction, there is a rational function $R_\zeta:\mathbb P^1\to \mathbb P^1$ of degree three such that $R_\zeta(f_\zeta^1(z))=z$.
Therefore, $f_\zeta^1$ is irrational, i.e., $f_\zeta^1(w)$ is not a rational function in $w\in \Delta\subset \mathbb C$.
Since any totally geodesic holomorphic isometry $f:(\Delta,g_\Delta)\to (\Delta,g_\Delta)\times (\mathbb B^2,g_{\mathbb B^2})$ with $f(0)={\bf 0}$ is the restriction of some linear map $\mathbb C\to \mathbb C^3$ (cf. Mok \cite[p.\;1646]{Mok12}), the constructed holomorphic isometry $f_\zeta$ is not totally geodesic.

Now, we show that $f_\zeta$ extends holomorphically to a neighborhood of $\overline\Delta$.
Note that the only possible singularities of $f_\zeta$ lying on the unit circle $\partial\Delta$ are branch points or poles of the component functions of $f_\zeta$. One observes from $R_\zeta(f_\zeta^1(w))=w$ and $f_\zeta^1:\Delta\to \Delta$ that $f_\zeta^1$ does not have any pole on the unit circle $\partial\Delta$.
For the rest of the proof, we follow the notations in the proof of Proposition \ref{Pro:Ex_1PFICHIndim}.
Note that a branch point of $f_\zeta^1$ is a branch point of the $3$-sheeted branched covering $R_\zeta:\mathbb P^1 \to \mathbb P^1$.
Since $0<|\zeta|<{1\over 3}$, the ramification points $a_+(\zeta)$ and $a_-(\zeta)$ of $R_\zeta$ satisfy $a_+(\zeta)={1\over \overline{a_-(\zeta)}}$ and $0<|a_+(\zeta)|<1$ (cf. the proof of Proposition \ref{Pro:Ex_1PFICHIndim}).
Moreover, $a_+(\zeta)$ and $a_-(\zeta)$ should lie in different fibers of $R_\zeta$ because $R_\zeta$ is $3$-sheeted so that the fiber of each branch point of $R_\zeta$ contains precisely one ramification point, i.e.,
$R_\zeta (a_-(\zeta))\neq R_\zeta(a_+(\zeta))$.
In addition, $a_-(\zeta)={1\over \overline{a_+(\zeta)}}$ so that $R_\zeta (a_-(\zeta)) = {1\over \overline{R_\zeta(a_+(\zeta))}}$ because $R_\zeta\left({1\over \overline z}\right) = {1\over \overline{R_\zeta(z)}}$.
If $R_\zeta (a_-(\zeta))\in \partial\Delta$, then we would have $R_\zeta (a_-(\zeta))=R_\zeta(a_+(\zeta))$, a plain contradiction.
Thus, we have $R_\zeta (a_+(\zeta)),\;R_\zeta (a_-(\zeta))\not\in \partial\Delta$.
In particular, all branch points of $f_\zeta^1$ are outside the closed unit disk $\overline\Delta$ so that $f_\zeta^1$ extends holomorphically to a neighborhood of $\overline\Delta$.

From the construction, we have $f_{\zeta}^{2,j}=R_j\circ f_\zeta^1$ for some rational function $R_j:\mathbb P^1 \to \mathbb P^1$, $j=1,2$.
Thus, it remains to show that $f_{\zeta}^{2,j}$ does not have any pole lying on the unit circle $\partial\Delta$ for $j=1,2$.
Actually, the set $P_j$ of all poles of $R_j$ is a subset of the set $P$ of all poles of $R_\zeta$ from the construction, $j=1,2$.
Then, we have $f_\zeta^1(b)\not\in P_j$ for any $b\in \partial\Delta$ and any $j$ because $R_\zeta(f_\zeta^1(b)) = b$ and $P_j\subset P$.
Thus, each $f_{\zeta}^{2,j}$ extends holomorphically to a neighborhood of $\overline\Delta$, $j=1,2$.
The proof is complete.
\end{proof}

\begin{rema}
\begin{enumerate}
\item[(1)]
This corollary shows that Theorem 3 in \cite[p.\;2637]{MN09} could not be generalized to the case where the target is product of complex unit balls $\mathbb B^{N_1}\times\cdots \times \mathbb B^{N_p}$, $p\ge 2$, with $N_j\ge 2$ for some $j$, $1\le j\le p$.
Nevertheless, one needs to impose the single valuedness of the holomorphic isometries in order to generalize Theorem 3 in \cite{MN09} (cf. Theorem \ref{Thm_ExtHolo}).
\item[(2)] As an application of the existence of the holomorphic isometry from $\Delta$ to $\Delta \times \mathbb B^2$ with irrational component function(s) that extends holomorphically to a neighborhood of $\overline\Delta$, Xiao and the two authors discovered that there exist proper holomorphic maps from $\Delta$ to $\mathbb B^2$ that are algebraic and extend holomorphically to a neighborhood of $\overline\Delta$ but are not rational. In fact, $f_\zeta^1$ maps $\overline\Delta$ onto a compact set $E \subset \Delta$ for any $0<|\zeta|<{1\over 3}$, where $f_\zeta^1: \Delta \to \Delta$ is the map constructed in the above corollary.
This is a new phenomenon in the rank 1 case (see \cite{CXY17} for more details).
\end{enumerate}
\end{rema}
\section{New examples of holomorphic isometries from the Poincar\'e disk into certain irreducible bounded symmetric domains of rank $\ge 2$}
In this section, we provide more applications of Proposition \ref{Pro:Ex_1PFICHIndim} and our study on holomorphic isometries from $(\Delta,g_\Delta)$ to $(\Delta,g_\Delta)\times(\mathbb B^n,g_{\mathbb B^n})$, $n\ge 2$,
to the study of holomorphic isometries from $(\Delta,g_\Delta)$ to $(\Omega,g_\Omega)$ for any irreducible bounded symmetric domain $\Omega$ of rank $\ge 2$.

Firstly, we construct examples of non-standard holomorphic isometries from the Poincar\'e disk into irreducible bounded symmetric domains of rank $\ge 2$ which are irrational and extend holomorphically to a neighborhood of the closed unit disk $\overline\Delta$.
It is well-known that any irreducible bounded symmetric domain is biholomorphic to one of the following:
\[ D^{\mathrm{I}}_{p,q} := \left\{Z\in M(p,q;\mathbb C) : {\bf I_q} - \overline{Z}^T Z > 0 \right\},\quad p, q \ge 1, \]
\[ D^{\mathrm{II}}_m := \left\{ Z\in D^{\mathrm{I}}_{m,m}: Z=-Z^T\right\},\quad m\ge 2, \]
\[ D^{\mathrm{III}}_m := \left\{ Z\in D^{\mathrm{I}}_{m,m}: Z=Z^T\right\},\quad m\ge 1, \]
\[ D^{\mathrm{IV}}_n:=\left\{(z_1,\ldots,z_n)\in \mathbb C^n:\sum_{j=1}^n |z_j|^2 < 2,\;\sum_{j=1}^n |z_j|^2 < 1 +\left|{1\over 2}\sum_{j=1}^n z_j^2 \right|^2\right\},\; n\ge 3, \]
\[ D^{\mathrm{V}}\cong E_6/SO(10)\cdot SO(2) \text{ of complex dimension } 16,\]
\[ D^{\mathrm{VI}}\cong E_7/E_6\cdot SO(2)  \text{ of complex dimension } 27 \]
(cf. \cite{Wo72, Mok89}).
The corresponding K\"ahler form $\omega_{g_{D^{\mathrm{IV}}_n}}$ of $(D^{\mathrm{IV}}_n,g_{D^{\mathrm{IV}}_n})$ is given by
\[ \omega_{g_{D^{\mathrm{IV}}_n}} = -\sqrt{-1}
\partial\overline \partial \log \left(1-\sum_{j=1}^n |z_j|^2 +\left|{1\over 2}\sum_{j=1}^n z_j^2 \right|^2\right). \]

Throughout this section, we denote by $f_{\zeta,2}=f_\zeta:(\Delta,g_\Delta)\to (\Delta,g_\Delta)\times (\mathbb B^{2},g_{\mathbb B^{2}})$ the non-standard holomorphic isometry obtained in Corollary \ref{Cor:ExHolo_n=2} for $0<|\zeta|<{1\over 3}$.
For $n\ge 3$, let $f_{\zeta,n}:(\Delta,g_\Delta)\to (\Delta,g_\Delta)\times (\mathbb B^{n},g_{\mathbb B^{n}})$ be the non-standard holomorphic isometry defined by $f_{\zeta,n}(w):=(f_{\zeta,2}(w),0,\ldots,0)$
for $0<|\zeta|<{1\over 3}$.
Then, the holomorphic isometry $f_{\zeta,n}$ extends holomorphically to a neighborhood of $\overline\Delta$ and has irrational component function(s) for any $n\ge 2$ and any $\zeta\in \mathbb C$ such that $0<|\zeta|<{1\over 3}$ by Corollary \ref{Cor:ExHolo_n=2}.

From now on, we fix $\zeta\in \mathbb C$ with $0<|\zeta|<{1\over 3}$.
Let $\Omega\Subset \mathbb C^N$ be an irreducible bounded symmetric domain of rank $r\ge 2$ in its Harish-Chandra realization.
Then, it follows from Wolf \cite{Wo72} that the rank-$1$ boundary component of $\Omega$ is the complex unit ball of complex dimension $n_{r-1}(\Omega)$, where $n_{r-1}(\Omega)$ is the $(r-1)$-th null dimension of $\Omega$ (cf. \cite{Mok89}).
Moreover, there exists a totally geodesic holomorphic isometric embedding $G:(\Delta,g_\Delta)\times (\mathbb B^{n_{r-1}(\Omega)},g_{\mathbb B^{n_{r-1}(\Omega)}})\hookrightarrow (\Omega,g_\Omega)$ (cf. \cite{Wo72}).
We may assume without loss of generality that $G({\bf 0})={\bf 0}$ and $G$ is the restriction of the linear map from $\mathbb C^{n_{r-1}(\Omega)+1}$ to $\mathbb C^N$.

If $\Omega$ is of non-tube type, then it follows from Mok \cite{Mok02} that $n:=n_{r-1}(\Omega)\ge 2$.
Therefore, the composed map $G\circ f_{\zeta,n}:(\Delta,g_\Delta)\to (\Omega,g_\Omega)$ is a non-standard holomorphic isometry which extends holomorphically to a neighborhood of $\overline\Delta$ and is irrational.

Now, we assume that $\Omega$ is of tube type.
Then, $\Omega$ is biholomorphic to either $D^{\mathrm{I}}_{p,p}$ for some $p\ge 2$, $D^{\mathrm{II}}_{2k}$ for some $k\ge 2$, $D^{\mathrm{III}}_m$ for some $m\ge 2$, $D^{\mathrm{IV}}_n$ for some $n\ge 3$ or $D^{\mathrm{VI}}$.
Note that there is a totally geodesic holomorphic isometric embedding $G_1:(\Delta,g_\Delta)\times (\Omega',g_{\Omega'})\hookrightarrow (\Omega,g_\Omega)$, where $\Omega'\subset \Omega$ is the maximal characteristic symmetric subdomain and is an irreducible bounded symmetric domain of rank $\mathrm{rank}(\Omega)-1=r-1$.
Moreover, it follows from Wolf \cite{Wo72} that
\[ \Omega'\cong \begin{cases}
D^{\mathrm{I}}_{p-1,p-1} & \text{if } \Omega\cong D^{\mathrm{I}}_{p,p},\; p\ge 2.\\
D^{\mathrm{II}}_{2k-2} & \text{if } \Omega\cong D^{\mathrm{II}}_{2k},\;k\ge 2.\\
D^{\mathrm{III}}_{m-1} & \text{if } \Omega\cong D^{\mathrm{III}}_{m},\;m\ge 2.\\
\Delta& \text{if } \Omega\cong D^{\mathrm{IV}}_{n},\;n\ge 3.\\
D^{\mathrm{IV}}_{10} & \text{if } \Omega\cong D^{\mathrm{VI}}.
\end{cases} \]
We observe that any irreducible bounded symmetric domain of rank $2$ and of tube type is isometrically biholomorphic to a type-$\mathrm{IV}$ domain $D^{\mathrm{IV}}_n$ for some $n\ge 3$.
Thus, we further restrict to the case where $\Omega$ is of rank $\ge 3$ because the case where $\Omega\cong D^{\mathrm{IV}}_n$ for some $n\ge 3$ will be done by another method.

If $\Omega$ is biholomorphic to either $D^{\mathrm{I}}_{p,p}$ for some $p\ge 3$, $D^{\mathrm{II}}_{2k}$ for some $k\ge 3$, or $D^{\mathrm{VI}}$, then it follows from \cite[p.\;1226]{HT02} that there exists a totally geodesic holomorphic isometric embedding $G_2:(\mathbb B^m,g_{\mathbb B^m})\hookrightarrow (\Omega',g_{\Omega'})$ for some $m\ge 2$.
Thus, there is a totally geodesic holomorphic isometric embedding $G_3:(\Delta,g_\Delta)\times (\mathbb B^m,g_{\mathbb B^m})\hookrightarrow (\Omega,g_\Omega)$ for some $m\ge 2$.
We may assume without loss of generality that $G_3({\bf 0})={\bf 0}$ and $G_3$ is the restriction of the linear map from $\mathbb C^{m+1}$ to $\mathbb C^N$.
Similar to the case of irreducible bounded symmetric domains of non-tube type, the composed map $G_3 \circ f_{\zeta,m}:(\Delta,g_\Delta)\to (\Omega,g_\Omega)$ is a non-standard holomorphic isometry which extends holomorphically to a neighborhood of $\overline\Delta$ and is irrational.

We now consider the case of $\Omega\cong D^{\mathrm{III}}_m$ for some $m\ge 4$.
Recall that there is a totally geodesic holomorphic isometric embedding $(\Delta,g_\Delta)\times (D^{\mathrm{III}}_{m-1},g_{D^{\mathrm{III}}_{m-1}})\hookrightarrow (D^{\mathrm{III}}_{m},g_{D^{\mathrm{III}}_{m}})$ given by $(w,Z)\mapsto \begin{bmatrix}
w & {\bf 0}\\
{\bf 0} & Z
\end{bmatrix}$.
On the other hand, Xiao-Yuan \cite[Theorem 6.13]{XY16b} constructed a polynomial holomorphic isometry $F_n:(\mathbb B^{n-1},g_{\mathbb B^{n-1}})\to (D^{\mathrm{III}}_{n},g_{D^{\mathrm{III}}_{n}})$ for $n\ge 2$.
Here $F$ being a polynomial holomorphic isometry means that $F$ is a holomorphic isometry and each component function of $F$ is a polynomial in $z\in \mathbb B^{n-1}\Subset \mathbb C^{n-1}$.
In particular, there is a polynomial holomorphic isometry $F^{\mathrm{III}}_m:(\Delta,g_\Delta)\times (\mathbb B^{m-2},g_{\mathbb B^{m-2}})\hookrightarrow (D^{\mathrm{III}}_{m},g_{D^{\mathrm{III}}_{m}})$ given by
\[ F^{\mathrm{III}}_m(w,z_1,\ldots,z_{m-2})
=\begin{bmatrix}
w & {\bf 0}\\
{\bf 0} & F_{m-1}(z_1,\ldots,z_{m-2})
\end{bmatrix} \]
for $m\ge 4$.
Then, the composed map $F^{\mathrm{III}}_m\circ f_{\zeta,m-2}:(\Delta,g_\Delta)\to (D^{\mathrm{III}}_{m},g_{D^{\mathrm{III}}_{m}})$, $m\ge 4$, is a non-standard holomorphic isometry which extends holomorphically to a neighborhood of $\overline\Delta$ and is irrational.

Now, we consider the case of $\Omega\cong D^{\mathrm{III}}_3$.
Note that there is a holomorphic isometry $\nu: (D^{\mathrm{IV}}_3,g_{D^{\mathrm{IV}}_3})\to (D^{\mathrm{III}}_2,g_{D^{\mathrm{III}}_2})$ which is a biholomorphism.
We may assume $\nu({\bf 0})={\bf 0}$ without loss of generality and thus $\nu$ is the restriction of a linear map $\mathbb C^3 \to M_s(2,2;\mathbb C):=\{Z\in M(2,2;\mathbb C): Z=Z^T\}\cong \mathbb C^3$ (cf. \cite[Proposition 3.1.1]{Mok12}).
On the other hand, Mok \cite{Mok16} constructed a non-standard holomorphic isometry $\widetilde F:(\mathbb B^2,g_{\mathbb B^2}) \to (D^{\mathrm{IV}}_3,g_{D^{\mathrm{IV}}_3})$.
Later on, Xiao-Yuan \cite[p.\;30]{XY16b} have also written down the map $\widetilde F$ explicitly and the only possible singularity of $\widetilde F$ at the boundary $\partial\mathbb B^2$ is the point $(0,1)$.
In particular, we obtain a holomorphic isometry
$\widetilde F^{\mathrm{III}}_3:(\Delta,g_\Delta)\times (\mathbb B^2,g_{\mathbb B^2})\to (D^{\mathrm{III}}_3,g_{D^{\mathrm{III}}_3})$ given by
\[ \widetilde F^{\mathrm{III}}_3(w,z_1,z_2)
= \begin{bmatrix}
w &{\bf 0}\\
{\bf 0}& \nu\circ \widetilde F(z_1,z_2)
\end{bmatrix} \]
We claim that the composed map $\widetilde F^{\mathrm{III}}_3\circ f_{\zeta,2}$ extends holomorphically to a neighborhood of $\overline\Delta$.
Actually, it suffices to show that $f_{\zeta,2}(w)$ does not tend to $(w_0,0,1)$ for any $w_0\in \partial\Delta$.
Write $f_{\zeta,2}(w)
:=(f_{\zeta,2}^1(w),g_{\zeta,2}^1(w),g_{\zeta,2}^2(w))$.
Then, we have
\[ g_{\zeta,2}^2(w) = {\sqrt{1-|\zeta|^2} f_{\zeta,2}^1(w)\over f_{\zeta,2}^1(w)-\zeta} \]
from the construction.
If there is $w_0\in \partial\Delta$ such that $g_{\zeta,2}^2(w)\to 1$ as $w\to w_0$, then we would have
$f_{\zeta,2}^1(w)\to {\zeta\over 1-\sqrt{1-|\zeta|^2}}$ as $w\to w_0$.
We write $f_{\zeta,2}^1(w_0):={\zeta\over 1-\sqrt{1-|\zeta|^2}}$.
Then, we have $|f_{\zeta,2}^1(w_0)|\le 1$ so that $|\zeta|\le 1-\sqrt{1-|\zeta|^2}$.
Since $0<|\zeta|<{1\over 3}$, $|\zeta|\le 1-\sqrt{1-|\zeta|^2}$ implies that
$\sqrt{1-|\zeta|^2} \le 1-|\zeta| \le 1-|\zeta|^2$, i.e.,
$1-\sqrt{1-|\zeta|^2}\le 0$, a plain contradiction.
Therefore, the composed map $\widetilde F^{\mathrm{III}}_3\circ f_{\zeta,2}:(\Delta,g_\Delta)\to (D^{\mathrm{III}}_3,g_{D^{\mathrm{III}}_3})$ is a non-standard holomorphic isometry which extends holomorphically to a neighborhood of $\overline\Delta$ and is irrational.

\medskip
It remains to consider the case of $\Omega\cong D^{\mathrm{IV}}_n$ for some $n\ge 3$.
From \cite{CM17, XY16a}, there is a non-standard holomorphic isometry
$F:(\mathbb B^2,g_{\mathbb B^2})\to (D^{\mathrm{IV}}_3,g_{D^{\mathrm{IV}}_3})$ given by
\[ F(z_1,z_2):=\left(z_1,z_2,1-\sqrt{1-\sum_{j=1}^2 z_j^2}\right). \]
Let $h:(\Delta,g_\Delta)\to (\mathbb B^2,g_{\mathbb B^2})$ be the holomorphic isometry given by $h(w)=\left({\sqrt{-1}\over 4} w, {\sqrt{15}w\over 4}\right)$.
Then, we have
$F\circ h(w) = \left({\sqrt{-1}\over 4} w, {\sqrt{15}\over 4} w
,1-\sqrt{1-{7\over 8} w^2}\right)$.
Note that $1-{7\over 8} w^2\neq 0$ for any $w\in \overline\Delta$.
Define the map $\widetilde F^{\mathrm{IV}}_n:\Delta\to D^{\mathrm{IV}}_n$, $n\ge 3$, by
\[ \widetilde F^{\mathrm{IV}}_n(w)
= \left({\sqrt{-1}\over 4} w, {\sqrt{15}\over 4}w
,1-\sqrt{1-{7\over 8} w^2},0,\ldots,0\right). \]
Then, $\widetilde F^{\mathrm{IV}}_n:(\Delta,g_\Delta)\to (D^{\mathrm{IV}}_n,g_{D^{\mathrm{IV}}_n})$ is a non-standard holomorphic isometry which extends holomorphically to a neighborhood of $\overline\Delta$ and is irrational.

In short, we have shown the following:
\begin{theo}\label{Pro:NewEx_ExHolo_irrat}
Let $\Omega\Subset \mathbb C^N$ be an irreducible bounded symmetric domain of rank $\ge 2$ in its Harish-Chandra realization.
Then, there is a non-standard holomorphic isometry $F:(\Delta,g_\Delta)\to (\Omega,g_\Omega)$ which extends holomorphically to a neighborhood of the closed unit disk $\overline\Delta$ and is irrational, i.e., some component functions of $F=F(w)$ are not rational functions of $w\in \Delta\subset \mathbb C$.
\end{theo}
\begin{rema}
\begin{enumerate}
\item[(1)] This actually answers Problem 5.2.2. in Mok \cite{Mok11} in the negative when $\Omega$ is an irreducible bounded symmetric domain of rank $\ge 2$ and the normalizing constant $\lambda$ (cf. \cite{Mok11}) is the minimal possible one.
In addition, our examples of holomorphic isometries from $(\Delta,g_\Delta)$ to $(\Omega,g_\Omega)$ are all irrational.
\item[(2)] In \cite{XY16b}, Xiao-Yuan constructed non-standard holomorphic isometries from $(\Delta,g_\Delta)$ to $(\Omega,g_\Omega)$ which extend holomorphically to a neighborhood of $\overline\Delta$ and are rational, where $\Omega$ is any classical irreducible bounded symmetric domain of rank $\ge 2$ (cf. \cite{Mok16} as well).
Nevertheless, Theorem \ref{Pro:NewEx_ExHolo_irrat} shows that for any irreducible bounded symmetric domain $\Omega$ of rank $\ge 2$, there are non-standard holomorphic isometries from $(\Delta,g_\Delta)$ to $(\Omega,g_\Omega)$ which extend holomorphically to a neighborhood of $\overline\Delta$ and are irrational.
\item[(3)] Suppose that $\Omega$ is an irreducible bounded symmetric domain of rank $2$ and $F$ is as in Theorem \ref{Pro:NewEx_ExHolo_irrat}.
Then, the image of $F$ is not contained in any totally geodesic bidisk $\Delta^2\cong\Pi\subset \Omega$.
Otherwise, $F$ would be the composition of the square root embedding from $\Delta$ to $\Delta^2$ and the totally geodesic holomorphic isometric embedding from $(\Delta,g_\Delta)\times(\Delta,g_\Delta)$ to $(\Omega,g_\Omega)$ so that $F$ has two distinct branch points on the unit circle $\partial\Delta$ (cf. \cite{Ng10}), which contradicts with the fact that $F$ extends holomorphically to a neighborhood of the closed unit disk $\overline\Delta$. This answers Problem 5.2.5 in \cite{Mok11} in the negative.
\end{enumerate}
\end{rema}

\begin{theo}\label{pro:FamHI_MuIC1_ClassicalDomain}
Let $\Omega\Subset \mathbb C^N$ be an irreducible bounded symmetric domain of classical type and of rank $\ge 2$ such that $\Omega\not\cong D^{\mathrm{IV}}_n$ for any $n\ge 3$.
Then, there is a family of holomorphic isometries
$F_{\zeta}:(\Delta,g_\Delta)\to (\Omega,g_\Omega)$, $\alpha<|\zeta|<1$ for some $\alpha\in (0,1)$, such that $F_\zeta$ and $F_{\zeta'}$ are incongruent to each other provided that $|\zeta|\neq |\zeta'|$.
In particular, there exists a real $1$-parameter family of mutually incongruent holomorphic isometries $\widetilde F_t:(\Delta,g_\Delta)\to (\Omega,g_\Omega)$, $t\in \mathbb R$.
\end{theo}
\begin{proof}
Let $\widetilde f_{\zeta,n}=f_{\zeta}:(\Delta,g_\Delta)\to (\Delta,g_\Delta)\times (\mathbb B^{n},g_{\mathbb B^{n}})$ be the holomorphic isometry constructed in the proof of Proposition \ref{Pro:Ex_1PFICHIndim} for ${{n}-1\over n+1}<|\zeta|<1$ and $n\ge 2$.
Under the assumption, $\Omega$ is either biholomorphic to (I) $D^{\mathrm{I}}_{p,q}$ for some $p$ and $q$ satisfying $q\ge p\ge 2$ and $(p,q)\neq (2,2)$, (II) $D^{\mathrm{II}}_m$ for some $m\ge 5$ or (III) $D^{\mathrm{III}}_m$ for some $m\ge 3$.
Then, we have constructed many holomorphic isometries
$F_{\zeta,m,\Omega}:= G\circ \widetilde f_{\zeta,m}:(\Delta,g_\Delta)\to (\Omega,g_\Omega)$, where $G:(\Delta,g_\Delta)\times (\mathbb B^m,g_{\mathbb B^m})\hookrightarrow (\Omega,g_\Omega)$ is a holomorphic isometric embedding for some $m\ge 2$.
Write $\widetilde f_{\zeta,n}:=(\widetilde f_{\zeta,n}^1,\widetilde f_{\zeta,n}^{2,1},\ldots,\widetilde f_{\zeta,n}^{2,n})$.

\paragraph{Case I: $\Omega=D^{\mathrm{I}}_{p,q}$, $q\ge p \ge 2$ and $(p,q)\neq (2,2)$.}
Let
\[ G(w,z_1,\ldots,z_{q-1}) = \begin{bmatrix}
w & 0 &\cdots & 0\\
0 & z_1 & \cdots & z_{q-1}\\
{\bf 0} & {\bf 0} &\cdots & {\bf 0}
\end{bmatrix}. \]
Then, we have
\[ F_{\zeta,q-1,D^{\mathrm{I}}_{p,q}}(w)=G\circ \widetilde f_{\zeta,q-1}(w)
= \begin{bmatrix}
\widetilde f_{\zeta,q-1}^1(w) & 0 &\cdots & 0\\
0 & \widetilde f_{\zeta,q-1}^{2,1}(w) & \cdots & \widetilde f_{\zeta,q-1}^{2,q-1}(w)\\
{\bf 0} & {\bf 0} &\cdots & {\bf 0}
\end{bmatrix}. \]
Note that automorphisms of $D^{\mathrm{I}}_{p,q}$ are of the form $Z\mapsto (AZ+B)(CZ+D)^{-1}$ for $\begin{bmatrix}
A&B\\
C&D
\end{bmatrix}\in SU(p,q)$.
Suppose that $F_{\zeta,q-1,D^{\mathrm{I}}_{p,q}}$ and $F_{\zeta',q-1,D^{\mathrm{I}}_{p,q}}$ are congruent to each other for some $\zeta,\zeta'\in \mathbb C$ such that
${q-2\over q}<|\zeta|,|\zeta'|<1$.
Then, we have
\[ \widetilde f_{\zeta,q-1}^1(\phi(w)) = {a \widetilde f_{\zeta',q-1}^1(w)+b\over c \widetilde f_{\zeta',q-1}^1(w) + d} \]
for some $a,b,c,d\in \mathbb C$ and some $\phi\in \mathrm{Aut}(\Delta)$.
But then since $\widetilde f_{\zeta,q-1}^1$ is a non-constant function, the matrix $\begin{pmatrix}
a&b\\c&d
\end{pmatrix}$ is invertible and $z\mapsto {a z+b\over c z + d}=:\psi(z)$ is an automorphism of $\mathbb P^1$.
Let $R_\xi:\mathbb P^1\to \mathbb P^1$ be the rational function such that $R_\xi(\widetilde f_{\xi,q-1}^1(w))=w$.
Note that $0$ is the only branch point of $R_\xi$ which lies in $\Delta$ for $\xi\in \mathbb C$ satisfying ${q-2\over q}<|\xi|<1$.
Then, we have
\begin{equation}\label{Eq_I_1}
\phi^{-1} \circ R_\zeta\circ \psi = R_{\zeta'}.
\end{equation}
Since $\phi^{-1}\in \mathrm{Aut}(\Delta)$ and it maps branch points of $R_\zeta$ to the branch points of $R_{\zeta'}$ by Eq. (\ref{Eq_I_1}), we have $\phi^{-1}(0)=0$ so that $\phi(z)=e^{i\theta}z$ for some $\theta\in [0,2\pi)$.
In particular, we have
\begin{equation}\label{Eq_I_2}
e^{-i\theta}R_\zeta(\psi(z)) \equiv R_{\zeta'}(z).
\end{equation}
Then, we have $R_\zeta(\psi(\infty))=\infty$ by Eq. (\ref{Eq_I_2}) so that $\psi(\infty)\in \{\infty,\zeta\}$.
Recall that $\zeta$ is a ramification point of $R_\zeta$.
Since $\psi$ maps the ramification points of $R_{\zeta'}$ to that of $R_\zeta$, we have $\psi(\infty)=\infty$.
Similarly, $R_\zeta(\psi(0)) = 0$ by Eq. (\ref{Eq_I_2}) so that $\psi(0)\in \left\{0,{1\over \overline\zeta}\right\}$.
Since ${1\over \overline\zeta}$ is a ramification point of $R_\zeta$, we have $\psi(0)=0$.
Therefore, $\psi(z)=az$ for some nonzero complex number $a$.
Moreover, we have $\psi(\zeta')=\zeta$ and $\psi\left({1\over \overline{\zeta'}}\right)={1\over \overline\zeta}$ by Eq. (\ref{Eq_I_2}) and comparing the set of zeros (resp. set of poles) of $R_\zeta$ with the set of zeros (resp. set of poles) of $R_{\zeta'}$. 
In particular, $a\zeta'=\zeta$ and $a{1\over \overline{\zeta'}} = {1\over \overline\zeta}$.
This implies that $|\zeta|=|\zeta'|$ and $|a|^2=1$.
Therefore, $F_{\zeta,q-1,D^{\mathrm{I}}_{p,q}}$ and $F_{\zeta',q-1,D^{\mathrm{I}}_{p,q}}$ are incongruent to each other for any $\zeta,\zeta'\in \mathbb C$ such that ${q-2\over q}<|\zeta|,|\zeta'|<1$ and $|\zeta|\neq |\zeta'|$.
The result follows in this case.

\paragraph{Case II: $\Omega=D^{\mathrm{II}}_{m}$, $m \ge 5$.}
We let
\[ G(w,z_1,\ldots,z_{m-3}) = \begin{bmatrix}
w\cdot J_1 & {\bf 0}\\
{\bf 0} & \widetilde G({\bf z})
\end{bmatrix}, \]
where $J_1:=\begin{pmatrix}
0 & 1\\
-1 & 0
\end{pmatrix}$ and $\widetilde G({\bf z})
:= \begin{bmatrix}
0 & {\bf z} \\
-{\bf z}^T & {\bf 0}
\end{bmatrix}$ with ${\bf z} = \begin{pmatrix}
z_1,\ldots,z_{m-3}
\end{pmatrix}$.
Then, we have
\[\begin{split}
 F_{\zeta,m-3,D^{\mathrm{II}}_{m}}(w)=&G\circ \widetilde f_{\zeta,m-3}(w)\\
=& \begin{bmatrix}
\widetilde f_{\zeta,m-3}^1(w) J_1 & {\bf 0}\\
{\bf 0} & \widetilde G \left(\widetilde f_{\zeta,m-3}^{2,1}(w),\ldots,\widetilde f_{\zeta,m-3}^{2,m-3}(w)\right)
\end{bmatrix}. \end{split}\]
Note that automorphisms of $D^{\mathrm{II}}_{m}$ are of the form $Z\mapsto (AZ+B)(CZ+D)^{-1}$ for $\begin{bmatrix}
A&B\\
C&D
\end{bmatrix}\in SO^*(2m)$.
By the same argument as in Case $\mathrm{I}$, we see that $F_{\zeta,m-3,D^{\mathrm{II}}_{m}}$ and $F_{\zeta',m-3,D^{\mathrm{II}}_{m}}$ are incongruent to each other for any $\zeta,\zeta'\in \mathbb C$ such that
${m-4\over m-2}<|\zeta|,|\zeta'|<1$ and $|\zeta|\neq |\zeta'|$,
where $m\ge 5$.
The result follows in this case.

\paragraph{Case III: $\Omega=D^{\mathrm{III}}_{m}$, $m \ge 3$.}
We let
\[ G(w,z_1,\ldots,z_{m-1}) = \begin{bmatrix}
w & {\bf 0}\\
{\bf 0} & \widetilde G(z)
\end{bmatrix}, \]
where $\widetilde G:(\mathbb B^{m-1},g_{\mathbb B^{m-1}})\to (D^{\mathrm{III}}_{m-1},g_{D^{\mathrm{III}}_{m-1}})$ is the holomorphic isometric embedding constructed by Mok \cite{Mok16} and $z=(z_1,\ldots,z_{m-1})$.
Then, we have
\[\begin{split}
 F_{\zeta,m-1,D^{\mathrm{III}}_{m}}(w)=&G\circ \widetilde f_{\zeta,m-1}(w)\\
=& \begin{bmatrix}
\widetilde f_{\zeta,m-1}^1(w) & {\bf 0}\\
{\bf 0} & \widetilde G \left(\widetilde f_{\zeta,m-1}^{2,1}(w),\ldots,\widetilde f_{\zeta,m-1}^{2,m-1}(w)\right)
\end{bmatrix}. \end{split}\]
Note that automorphisms of $D^{\mathrm{III}}_{m}$ are of the form $Z\mapsto (AZ+B)(CZ+D)^{-1}$ for $\begin{bmatrix}
A&B\\
C&D
\end{bmatrix}\in Sp(m,\mathbb R)$.
By the same argument as in Case $\mathrm{I}$, we see that $F_{\zeta,m-1,D^{\mathrm{III}}_{m}}$ and $F_{\zeta',m-1,D^{\mathrm{III}}_{m}}$ are incongruent to each other for any $\zeta,\zeta'\in \mathbb C$ such that ${m-2\over m}<|\zeta|,|\zeta'|<1$ and $|\zeta|\neq |\zeta'|$, where $m\ge 3$.
The result follows in this case.
\end{proof}
\begin{rema}\label{remark3.24}
\begin{enumerate}
\item[(1)]
When $\Omega$ is any irreducible bounded symmetric domain of rank $\ge 3$, it follows from Mok \cite[Theorem 3.2.1]{Mok12} that there is a real $1$-parameter family of mutually incongruent holomorphic isometries $H_t:(\Delta,g_\Delta)\to (\Omega,g_{\Omega})$, $t\in \mathbb R$, by using the Polydisk Theorem and the arguments in the proof of \cite[Theorem 3.2.1]{Mok12}.
Now, we restrict to the case where the target $\Omega$ is of classical type and of rank $\ge 3$.
Then, we observe that the holomorphic isometry $\widetilde F_t$ has at most two distinct branch points on the unit circle $\partial\Delta$ while $H_t$ has four distinct branch points on $\partial\Delta$ (cf. the proof of \cite[Theorem 3.2.1]{Mok12}).
This shows that $\widetilde F_t$ is incongruent to $H_{t'}$ for any  $t,t'\in \mathbb R$.
\item[(2)] In \cite{XY16b}, Xiao-Yuan proved the existence of a real $1$-parameter family of mutually incongruent holomorphic isometries $\widetilde F_t:(\Delta,g_\Delta)$ $\to$ $(D^{\mathrm{IV}}_n,g_{D^{\mathrm{IV}}_n})$, $t\in \mathbb R$.
On the other hand, our result in Theorem \ref{pro:FamHI_MuIC1_ClassicalDomain} yields a new phenomenon when the target $\Omega$ is biholomorphic to $D^{\mathrm{I}}_{2,q}$ or $D^{\mathrm{II}}_{5}$, where $q\ge 3$ is any integer.
\end{enumerate}
\end{rema}

\section{Rigidity of rational holomorphic isometries from the unit disk into a product of complex unit balls}
Recall that we have an example of non-standard holomorphic isometry $(\Delta,g_\Delta)\to (\Delta,g_\Delta)\times (\mathbb B^n,g_{\mathbb B^n})$, $n\ge 2$, which extends holomorphically to a neighborhood of the closed unit disk $\overline\Delta$.
In particular, one has to impose some stronger assumptions in order to generalize Theorem 3 in \cite[p.\;2637]{MN09} and obtain the rigidity of a certain class of holomorphic isometries from the unit disk into a product of complex unit balls.
\begin{theo}\label{Thm_ExtHolo}
Let $f=(f_1,\ldots,f_m):(\Delta,g_\Delta) \to (\mathbb B^{N_1},\lambda_1g_{\mathbb B^{N_1}}) \times \cdots \times (\mathbb B^{N_m},\lambda_mg_{\mathbb B^{N_m}})$ be a holomorphic isometry
such that $f_j$ is a non-constant map for $1\le j\le m$, where $\lambda_j$, $1\le j\le m$, are positive real constants. 
If $f$ is rational, i.e., each component function of $f$ is a rational function in $w\in \Delta\subset \mathbb C$, then $f_j:(\Delta,g_\Delta)\to (\mathbb B^{N_j},g_{\mathbb B^{N_j}})$ is a (totally geodesic) holomorphic isometry for $1\le j\le m$ and $\sum_{j=1}^m\lambda_j=1$ so that $f$ is totally geodesic.
\end{theo}
\begin{proof}
Assume without loss of generality that $f(0)={\bf 0}$.
Then, we have the functional equation
\[ \prod_{j=1}^m (1-\lVert f_j(w) \rVert^2)^{\lambda_j}
= 1-|w|^2 \]
for any $w\in \Delta$.
Write $f_j:=(f_j^1,\ldots,f_j^{N_j})$ for $1\le j\le m$.
We claim that $f$ extends holomorphically to a neighborhood of $\overline\Delta$.
If there is $z_0\in \partial\Delta$ such that $z_0$ is a pole of $f_j^l$ for some $j$ and some $l$, $1\le l\le N_j$, then $1-\lVert f_j(w) \rVert^2$ tends to negative infinity as $w\to z_0$.
Take a simple continuous path $\gamma(t)$, $0\le t\le 1$, from $0$ to $z_0$ in $\Delta \cup \{z_0\}$ such that $\gamma(0)=0$ and $\gamma(1)=z_0$. We see that there is $t_0\in (0,1)$ such that $\lVert f_j(\gamma(t_0)) \rVert^2 = 1$ because $\lVert f_j(\gamma(0)) \rVert^2 = 0$, $\lVert f_j(\gamma(t)) \rVert^2$ is continuous on $[0,1)$ and $\lVert f_j(\gamma(t)) \rVert^2\to +\infty$ as $t\to 1$.
Note that $\gamma(t_0) \in \Delta$ so that $f$ is holomorphic around $\gamma(t_0)$.
This would lead to a contradiction by the functional equation since $1-|\gamma(t_0)|^2>0$.
In particular, the poles of the component functions of $f$ do not lie on $\partial\Delta$ so that $f$ extends holomorphically to a neighborhood of $\overline\Delta$.

If $f_j$ maps some point $w_0\in \partial\Delta$ to $\partial \mathbb B^{N_j}$, then $f_j$ actually maps some open arc $A\subset \partial\Delta$ containing $w_0$ into $\partial \mathbb B^{N_j}$ (cf. Proof of Theorem 2 in \cite[p.\;894]{Mok09}).
Since $f_j$ extends holomorphically to a neighborhood of $\overline\Delta$, $1-\lVert f_j(w) \rVert^2$ is a real-analytic function in a neighborhood of $\overline\Delta$.
Then, we have $1-\lVert f_j(w) \rVert^2\equiv 0$ on $\partial\Delta$ by the Identity Theorem for real-analytic functions.
In particular, $f_j:\Delta\to \mathbb B^{N_j}$ is a proper holomorphic map.

Assume without loss of generality that $f_j$ is proper for $1\le j\le k$ and $f_l(\overline\Delta)\subset \mathbb B^{N_l}$ for $k+1\le l\le m$ if $k<m$, where $k$ is some integer satisfying $1\le k\le m$.
We claim that $f_l$ is a constant map for $k+1\le l\le m$ if $k<m$.
Note that the functional equation
\[ \prod_{j=1}^m (1-\lVert f_j(w) \rVert^2)^{\lambda_j}
= 1-|w|^2 \]
holds true in $\mathbb{C}$ away from the poles of $f$ by the Identity Theorem for real-analytic functions.
Note that $1-|w|^2<0$ for $|w|^2>1$ so that $\prod_{j=1}^m (1-\lVert f_j(w) \rVert^2)^{\lambda_j}$ does not vanish.
If there is $w\in \mathbb C\smallsetminus\overline\Delta \subset \mathbb P^1$ such that $1-\lVert f_l(w) \rVert^2=0$ for some $l$, $k+1\le l\le m$, then $\prod_{1\le j\le m,\;j\neq l} (1-\lVert f_j(w) \rVert^2)^{\lambda_j}$ has a pole.
But then there are only finitely many poles among all $f_j$, $1\le j\le m$ with $j\neq l$.
By avoiding the poles of all $f_j$'s, we have $1-\lVert f_l(w) \rVert^2>0$ on a dense open subset of $\mathbb P^1$ for $k+1\le l\le m$.
This shows that $f_l(\mathbb P^1)\subset \overline{\mathbb B^{N_l}}$, which contradicts with the Liouville's Theorem unless $f_l$ is a constant map.
In particular, we would have $f_l\equiv {\bf 0}$ for $k+1\le l\le m$.
But then this contradicts with the assumption that $f_j$ is a non-constant map for $1\le j\le m$.
Therefore, we have $k=m$ so that $f_j:\Delta\to \mathbb B^{N_j}$ is a proper holomorphic map for $1\le j\le m$.

Now, we follow the arguments in \cite{YZ12}.
Then, ${1-\lVert f_j(w)\rVert^2\over 1-|w|^2}$ is smooth and nonzero in some neighborhood of $\overline\Delta$ for $1\le j\le m$.
From the functional equation, we have
\[ \left(1-\sum_{j=1}^m\lambda_j\right)
\sqrt{-1}\partial\overline\partial\log (1-|w|^2)
=\sum_{j=1}^m \lambda_j \sqrt{-1}\partial\overline\partial\log\left({1-\lVert f_j(w)\rVert^2\over 1-|w|^2}\right). \]
This implies that $1-\sum_{j=1}^m\lambda_j=0$ because $\partial\overline\partial\log (1-|w|^2)$ is singular around any point $b\in \partial\Delta$.
In addition, we have $\sqrt{-1}\partial\overline\partial\log\left({1-\lVert f_j(w)\rVert^2\over 1-|w|^2}\right)\ge 0$ by the Ahlfors-Schwarz Lemma so that
$\partial\overline\partial\log(1-\lVert f_j(w)\rVert^2)=\partial\overline\partial\log(1-|w|^2)$ for $1\le j\le m$.
Hence, $f_j:(\Delta,g_\Delta) \to (\mathbb B^{N_j},g_{\mathbb B^{N_j}})$ is a (totally geodesic) holomorphic isometry for $1\le j\le m$ and the proof is complete.
\end{proof}

\noindent Shan Tai Chan\\
Department of Mathematics, The University of Hong Kong, Pokfulam Road, Hong Kong\\
mastchan@hku.hk

\medskip\noindent
Yuan Yuan\\
Department of Mathematics, Syracuse University, Syracuse, NY 13244-1150, USA\\
yyuan05@syr.edu

\end{document}